\documentclass[11pt,reqno]{amsart}
\usepackage[margin=1in]{geometry}
\usepackage{algorithm,algpseudocode}
\usepackage[fleqn,tbtags]{mathtools}
\usepackage[shortlabels]{enumitem}
\usepackage{verbatim}
\usepackage{tikz}
\usepackage{tikz-cd}
\usepackage{filecontents}
\usepackage{placeins}
\usepackage{mlmodern}
\usepackage{extarrows}
\usepackage{graphicx}
\usepackage{calc}
\usepackage[labelformat=simple]{subcaption}
\usepackage{adjustbox}
\usepackage{amssymb}
\usepackage{eucal,bbm}
\usepackage[scr]{rsfso}
\usepackage{makecell}
\usepackage{tabu}
\usepackage{braket}
\usepackage{multirow}
\usepackage{rotating}
\usepackage{pdflscape}
\usepackage{pifont}
\usepackage{circledsteps}
\usepackage{hyperref,xcolor}
\usepackage{empheq}
\hypersetup{
    colorlinks,
    linkcolor={red!50!black},
    citecolor={blue!50!black},
    urlcolor={blue!80!black}
}

\usetikzlibrary{arrows.meta, positioning}
\let\oh=\circ
\newcommand{\ccirc}{\mathbin{\mathchoice
  {\xcirc\scriptstyle}
  {\xcirc\scriptstyle}
  {\xcirc\scriptscriptstyle}
  {\xcirc\scriptscriptstyle}
}}
\newcommand{\xcirc}[1]{\vcenter{\hbox{$#1\oh$}}}
\let\circ\ccirc

\newcommand\restr[2]{{
  \left.\kern-\nulldelimiterspace 
  #1 
  \vphantom{\big|} 
  \right|_{#2} 
  }}

\DeclareMathOperator{\diag}{diag}

\DeclareMathOperator{\rank}{rank}

\DeclareMathOperator{\tr}{tr}

\DeclareMathOperator{\Alt}{\mathsf{\Lambda}}
\DeclareMathOperator{\Sym}{\mathsf{S}}
\DeclareMathOperator{\Herm}{\mathsf{H}}

\DeclareMathOperator*{\argmin}{argmin}

\theoremstyle{definition}
\newtheorem{theorem}{Theorem}[section]
\newtheorem{definition}[theorem]{Definition}
\newtheorem{lemma}[theorem]{Lemma}
\newtheorem{corollary}[theorem]{Corollary}
\newtheorem{proposition}[theorem]{Proposition}

\numberwithin{equation}{section}

\newcommand{\tp}{{\scriptscriptstyle\mathsf{T}}}

\renewcommand{\a}{\mathfrak{a}}
\newcommand{\C}{\mathbb{C}}
\newcommand{\F}{\mathbb{F}}

\newcommand{\R}{\mathbb{R}}

\newcommand{\V}{\mathbb{V}}

\newcommand{\U}{\operatorname{U}}

\newcommand{\g}{\mathfrak{g}}
\renewcommand{\k}{\mathfrak{k}}
\renewcommand{\O}{\operatorname{O}}

\newcommand{\p}{\mathfrak{p}}
\newcommand{\G}{\mathsf{G}}
\newcommand{\K}{\mathsf{K}}

\newcommand{\ii}{{\ensuremath{\mathrm{i}}}}

\setlist[description]{font=\normalfont\itshape,labelindent=2em}

\allowdisplaybreaks

\title{Matrix nearness problems via Lie algebra}

\author{Rongbiao Thomas Wang}
\address{Computational and Applied Mathematics, University of Chicago, Chicago, IL 60637}
\email{rbwang@uchicago.edu}

\begin{document}

\begin{abstract}
As exemplified by the Eckart--Young theorem, matrix nearness problems with geometrically complicated
constraints can sometimes be solved by manipulating only singular values. In this article, we answer the question of when matrix nearness problems $D(A,X)$ or $D(X,A)$ can be solved by manipulating only the eigenvalues or singular values through the notion of isotropic divergences. We show under an isotropic divergence, a matrix nearness problem whose constraints respect the Lie group symmetry can always be reduced to an approximation problem on the eigenvalues or singular values of $A$, providing a systematic way to construct minimizers. We establish sufficient conditions of isotropy in the classes of Bregman, Jensen, and Calabi's diastasis and verify that numerous matrix divergences arising in applications are naturally isotropic, including the Frobenius, von Neumann, logdet, $S$-, quantum Jensen--Shannon, Bures--Wasserstein, quantum Hellinger, and quantum $q$-divergence, among others. The resulting framework opens the door for us to construct solutions to a myriad of matrix nearness problems including constraints on rank, spectrum, condition number, trace, determinant, as well as other physically motivated constraints.
\end{abstract}



\maketitle

\section{Introduction}

Given a matrix $A\in\Omega$, a divergence $D$ on $\Omega$, and a constraint set $\mathcal{S}\subseteq\Omega$, the goal of this article is to answer the question of when a minimizer to the matrix nearness problem
\begin{equation}\label{eq:near}
    \min_{X\in\mathcal{S}} D(A,X)
    \quad\text{or}\quad
    \min_{X\in\mathcal{S}} D(X,A),
\end{equation}
can be obtained by altering only the eigenvalues or singular values of $A$, a folklore strategy underlying classical solutions under various constraints, norms, and divergences \cite{Eck,FanHoffman1955,mir,Smolin2012,GVL1980,Higham1,Higham2}.

The mechanism is best illustrated via the Eckart--Young theorem \cite{Eck}, which solves the rank-constrained approximation problem
\begin{equation}\label{eq:rank}
    \min_{\rank(X)\leq r}\lVert A-X\rVert_{\mathsf F}, \qquad A \in \C^{m \times n}
\end{equation}
by truncating the singular values of $A$, despite the nonconvexity and complicated geometry of the feasible set \cite{Determinantal}. Let $A=U_A\Sigma_A V_A^*$ and $X=U_X\Sigma_X V_X^*$, be the singular values decomposition of $A$ and $X$. The key to the success of Eckart--Young theorem is that for all $X \in \C^{m \times n}$,
\begin{equation}\label{eq:ing2}
    \lVert A- U_A \Sigma_X V_A^*\rVert_{\mathsf F}
    =\lVert\sigma(A)-\sigma(X)\rVert_2
    \leq\lVert A-X\rVert_{\mathsf F},
\end{equation}
where $\sigma(A)$ and $\sigma(X)$ are respectively the vectors of singular values arranged in nonincreasing order for $A$ and $X$. Using \eqref{eq:ing2} and the fact that rank is invariant under left and right multiplications of unitary matrices, \eqref{eq:rank} is reduced to approximating $\sigma(A)$, i.e., 
\begin{equation}\label{eq:rank-equiv}
    \min_{\lvert \operatorname{supp}(x)\rvert \leq r} \lVert \sigma(A) - x \rVert_2, \qquad \operatorname{supp}(x) = \{i \in \{1,\dots, \min(m,n)\}: x_i \neq 0\}.
\end{equation}

Notice that \eqref{eq:ing2} is a property of the Frobenius norm independent of the constraint. Therefore, whenever an inequality of for the form \eqref{eq:ing2} holds for a divergence with respect to a matrix diagonalization, we may reduce a matrix nearness problem to a vector problem on the diagonalized matrix, whenever the feasible set respects the symmetry in the diagonalization. We formalize this reduction using Lie theory, following the long tradition of applying them to optimization and matrix theory \cite{Lewis1996,lewis2000convex,Tam1999,Schur,Horn,Kostant,Atiyah}. 

\subsection{Isotropic reduction}

In practice, we often want to approximate matrices from more structured domains such as Hermitian, symmetric, skew-symmetric, positive definite, and contraction matrices. Depending on the matrix space, the appropriate canonical form may arise from an eigenvalue, Autonne--Takagi \cite{Autonne15,Takagi1924}, or Youla \cite{MW1931,Youla} decompositions rather than the singular value decomposition. To treat these decompositions uniformly, we consider the \emph{isotropy representations} induced by symmetric spaces: a real vector spaces $\V$ of matrices over $\R$ or $\C$, equipped with an action of a compact Lie group $\K$, such as $\U(n)$, $\O(n)$, or their products. In this case, every $A\in\V$ admits an \emph{isotropy decomposition} $A=U\cdot A_+$, where $U\in\K$ and $A_+$ is a diagonal or block-diagonal matrix whose entries are the eigenvalues or singular values of $A$ \cite{ten,Tam1999,Knapp2002}. We write $\rho(A)$ for the corresponding vector of eigenvalues or singular values.

An $\K$-\emph{isotropic divergence} $D$ on an open $\K$-invariant domain $\Omega\subseteq\V$ is a $\K$-invariant divergence satisfying
\begin{equation}\label{eq:intro-alignment}
    D(A_+,B_+)\leq D(A,B)
    \qquad\text{for all }A,B\in\Omega,
\end{equation}
the generalization of \eqref{eq:ing2}, which gives the following reduction.

\begin{theorem}[Isotropic reduction, informally]
\label{thm:intro-reduction}
Let $D$ be an isotropic divergence on $\Omega$. If $\mathcal{S} \subseteq \Omega$ is $\K$-invariant, then for all $A \in \Omega$,
\begin{equation}\label{eq:near-reduced}
    \inf_{X\in\mathcal{S}}D(A,X)
    =\inf_{x\in\rho(\mathcal{S})}
        \widehat D\bigl(\rho(A),x\bigr), \qquad
    \inf_{X\in\mathcal{S}}D(X,A)
    =\inf_{x\in\rho(\mathcal{S})}
        \widehat D\bigl(x,\rho(A)\bigr),
\end{equation}
where $\widehat D\bigl(\rho(A),\rho(B)\bigr) \coloneqq D(A_+,B_+)$, is the \emph{reduced divergence} on vectors. Moreover, if $A = U \cdot A_+$ is the isotropic decomposition of $A$ and $x_*$ attains either infimum on the right-hand side of \eqref{eq:near-reduced}, then $X_* = U \cdot X_+$, where $X_+$ is a diagonal or block-diagonal matrix who entries are entries of $x_*$, attains the corresponding infimum on the left-hand side.
\end{theorem}

Therefore, whenever $D$ is isotropic and $\mathcal{S}$ is $\K$-invariant, an minimizer to \eqref{eq:near} can be obtained by only manipulating the eigenvalues or singular values of $A$.

\subsection{Isotropic divergences}

A variety of divergences arising in practice are isotropic. In Section~\ref{sec:divergences}, we establish sufficient conditions for isotropy in the classes of Bregman, Jensen divergences and Calabi diastases and verify the isotropy of a list of divergences arising from different applications and compute their reduced divergences. The list is summarized in Table~\ref{tab:isotropic-divergences}. 

\begin{table}[htb]
\centering
\footnotesize
\setlength{\tabcolsep}{2.1pt}
\tabulinesep=0.25ex
\renewcommand{\arraystretch}{1.05}
\resizebox{\textwidth}{!}{%
\begin{tabu}{@{}c|c|c|c|c|c@{}}
class & divergence & notation & $\Omega$
& $D(X,Y)$ & $\widehat D(u,v)$\\
\hline

\multirow{10}{*}{Bregman}
& Frobenius & $D_{\mathsf F}$ & $\p$
& $\tfrac12\|X-Y\|_{\mathsf F}^2$
& $\tfrac{c_\p}{2}\|u-v\|_2^2$\\

& von Neumann & $D_{\mathsf{vN}}$ & $\Herm_+(\F^n)$
& \makecell[c]{$\tr[X(\log X-\log Y)$\\
    $-X+Y]$}
& \makecell[c]{$\sum_i[u_i\log\frac{u_i}{v_i}$\\
    $-u_i+v_i]$}\\

& logdet & $D_{\mathsf{LD}}$ & $\Herm_+(\F^n)$
& \makecell[c]{$\tr(XY^{-1})-\log\det(XY^{-1})$\\
    $-n$}
& \makecell[c]{$\sum_i[\frac{u_i}{v_i}-\log\frac{u_i}{v_i}$\\
    $-1]$}\\

& & $D_{\mathsf{GKL}}$
& $\R^{m\times n}_{<1}$
& \makecell[c]{$\frac{1}{2}[\log\det(I-Y^*Y)$\\
 $- \log \det(I-X^*X)]$\\$
\Re\tr[(I-Y^*Y)^{-1}Y^*(Y-X)]$}
& \makecell[c]{$\tfrac12\sum_i[\log\frac{1-v_i^2}{1-u_i^2}$\\
    $-\frac{2v_i(u_i-v_i)}{1-v_i^2}]$}\\

& \makecell[c]{spectral\\hypentropy}
& $D_\beta$ & $\Herm(\F^n)$
& \makecell[c]{$D_{\varphi_\beta}(X,Y)$\\
    $\varphi_\beta(X)=\sum_i f_\beta(\lambda_i(X))$}
& \makecell[c]{$\sum_i[f_\beta(u_i)-f_\beta(v_i)$\\
    $-(u_i-v_i)\operatorname{arsinh}(v_i/\beta)]$}\\

\hline
\multirow{2}{*}{Jensen}
& S-divergence & $D_{\mathsf S}$ & $\Herm_+(\F^n)$
& \makecell[c]{$\log\det\frac{X+Y}{2}$\\
    $-\tfrac12\log\det(XY)$}
& $\sum_i\log\frac{u_i+v_i}{2\sqrt{u_iv_i}}$\\

& \makecell[c]{quantum\\Jensen--Shannon}
& $D_{\mathsf{QJS}}$ & $\Herm_+(\F^n)$
& \makecell[c]{$\tfrac12\tr(X\log X+Y\log Y)$\\
    $-\tr(M\log M)$}
& \makecell[c]{$\sum_i[\tfrac12(u_i\log u_i+v_i\log v_i)$\\
    $-m_i\log m_i]$}\\

\hline
\multirow{3}{*}{\makecell[c]{Calabi's\\diastasis}}
& & $C_f$
& \makecell[c]{$\Omega\subseteq\p_{<\mu}$\\
    $\p=\C^{m\times n},\ \Sym(\C^n),$\\
    $\Alt(\C^{2n})$}
& \makecell[c]{$\tr f(X^*X)+\tr f(Y^*Y)$\\
    $-2\Re\tr f(Y^*X)$}
& $c_\p\sum_{i,k\ge1}a_k(u_i^k-v_i^k)^2$\\

& & $C_{\mathsf B}$ & $\Sym(\C^n)_{<1}$
& \makecell[c]{$\log|\det(I-X^*Y)|$\\
    $-\tfrac12\log\det(I-X^*X)$\\
    $-\tfrac12\log\det(I-Y^*Y)$}
& \makecell[c]{$\sum_i\log\cosh(\operatorname{artanh}u_i$\\
    $-\operatorname{artanh}v_i)$}\\

& & $C_{\mathsf F}$ & $\Alt(\C^{2n})_{<1}$
& \makecell[c]{$\tfrac12\log\det(I+X^*X)$\\
    $+\tfrac12\log\det(I+Y^*Y)$\\
    $-\log|\det(I+X^*Y)|$}
& \makecell[c]{$-2\sum_i\log\cos(\arctan u_i$\\
    $-\arctan v_i)$}\\

\hline
\multirow{3}{*}{other}
& Bures--Wasserstein & $D_{\mathsf{BW}}$
& $\Herm_+(\F^n)$
& \makecell[c]{$\tr X+\tr Y$\\
    $-2\tr(X^{1/2}YX^{1/2})^{1/2}$}
& $\sum_i(\sqrt{u_i}-\sqrt{v_i})^2$\\

& quantum Hellinger & $D_{\mathsf{QH}}$
& $\Herm_+(\F^n)$
& $\lVert X^{1/2}-Y^{1/2}\rVert_{\mathsf F}^2$
& $\sum_i(\sqrt{u_i}-\sqrt{v_i})^2$\\

& quantum $q$-divergence & $D_{\mathsf{Q},q}$
& $\Herm_+(\F^n)$
& \makecell[c]{$\tfrac1{1-q}\tr[qX+(1-q)Y$\\
    $-X^qY^{1-q}]$}
& \makecell[c]{$\tfrac1{1-q}\sum_i[qu_i+(1-q)v_i$\\
    $-u_i^qv_i^{1-q}]$}
\end{tabu}
}
\caption{isotropic divergences and their reductions.}
\label{tab:isotropic-divergences}
\end{table}

Besides well-known divergences such as Frobenius, von Neumann, logdet, S-, quantum Jensen--Shannon, Bures--Wasserstein divergences, several surprising divergences in Table~\ref{tab:isotropic-divergences} arise from concrete statistical and physical models. The divergence $D_{\mathsf{GKL}}$ defined on real strict contractions arises from Gaussian models associated with canonical correlation analysis \cite{Hotelling1936,Anderson2003,BachJordan2005}. Calabi's diastasis, originally a purely geometric object from K\"ahler geometry \cite{Calabi1953}, measures the infidelity $C_{\mathsf B}$ and $C_{\mathsf F}$ for bosonic squeezed states and fermionic pairing states, respectively.

\subsection{Isotropic nearness problems}

To illustrate Theorem~\ref{thm:intro-reduction} in action, we construct in Section~\ref{sec:solve} the solutions to a range of matrix nearness problems over the divergences in Table~\ref{tab:isotropic-divergences} with $\mathcal{S}$ defined by
\begin{enumerate}[\normalfont(i)]
    \item rank constraint: $\rank(X)\leq r$ for a given integer $r>0$;
    \item spectral constraint:
    $\mathsf{\Lambda}(X)\cap[a,b]\neq\varnothing$ for given $a, b \in \R$;
    \item condition number constraint:
    $\sigma_1(X)\leq\kappa\sigma_n(X)$ for given $\kappa \geq 1$;
    \item trace constraint:
    $X\succeq0$, $\tr X=\tau$ for given $\tau > 0$;
    \item determinant constraint:
    $\det X=\delta$ for given $\delta > 0$;
\end{enumerate}
as well as constraints that arise from applications in quantum physics
\begin{enumerate}[\normalfont(i), start=6]
    \item bosonic squeezing resources:
    $\rank(X)\leq k$,
    $\tr[X^*X(I-X^*X)^{-1}]\leq N$ for given integers $0 < k \leq n$ and $N>0$;
    \item fermionic pairing resources:
    $\rank(X)\leq2k$,
    $\tr[X^*X(I+X^*X)^{-1}]\leq N$ for given integers $0 < k \leq n$ and $N>0$.
\end{enumerate}

To the best of our knowledge, most of the solutions obtained here have not previously appeared. We will point out the previously known cases in Section~\ref{sec:solve} throughout the discussions. We emphasize that we do not treat \eqref{eq:near} as a constrained optimization problem and then apply general purpose nonlinear, convex, or Riemannian optimization methods. While it is stated as an optimization problem, our solutions are constructive and rooted in numerical linear algebra and Lie theory, and relies decisively on the matrix structure in these problems. Tools from optimization theory are only used in solving the reduced problem \eqref{eq:near-reduced} analytically instead of algorithmically.

\subsection*{Notations and conventions} In this article, we focus on matrices taking values in the fields $\F= \R$ or $\C$, while the general theory extends to the case of quaternions. Let
\begin{align*}
    \mathsf{S}(\F^n) &\coloneqq \{ A \in \F^{n \times n} : A = A^\tp \}, \\
    \mathsf{\Lambda}(\F^n) &\coloneqq \{ A \in \F^{n \times n} : A = -A^\tp \},  \\
    \mathsf{H}(\F^n) &\coloneqq \{ A \in \F^{n \times n} : A = A^* \}
\end{align*}
denote the vector spaces of symmetric, skew-symmetric, and Hermitian matrices over $\F$ respectively. Notice that $\Herm(\R^n) = \Sym(\R^n)$ while $\Herm(\C^n) \neq \Sym(\C^n)$. We denote the set of positive definite matrices over $\F$ by $\Herm_+ (\F^n)$.

For $X \in \Herm(\F^n)$, we denote its spectrum by $\mathsf{\Lambda}(X) = \{\lambda_1(X),\dots,\lambda_n(X)\} \subseteq \R$, where eigenvalues are indexed in nonincreasing order and $\lambda(X) \coloneqq \begin{bmatrix}
    \lambda_1(X) &\cdots &\lambda_n(X)
\end{bmatrix}^\tp$. Similarly, for $Y \in \F^{m \times n}$, we denotes the set of its singular values by $\mathsf{\Sigma}(X) = \{\sigma_1(X),\dots,\sigma_n(X)\} \subseteq \R_{\geq 0}$, where singular values are arranged in nonincreasing order and $\sigma(X) \coloneqq \begin{bmatrix}
    \sigma_1(X) &\cdots &\sigma_n(X)
\end{bmatrix}^\tp$.

For a space of rectangle matrices $\F^{m \times n}$, we assume without loss of generality that $m \geq n$. In this case, we denote by $\diag_{m,n}(v)$ the matrix $\bigl[\begin{smallmatrix} \diag(v) \\0_{m-n \times n} \end{smallmatrix}\bigr]\in \F^{m \times n}$. Let $J = \left[\begin{smallmatrix} 0 & 1 \\ -1 & 0 \end{smallmatrix}\right]$.

\section{Lie-theoretic background}

A Riemannian symmetric space is a homogeneous space $\G/\K$ if there exists an involutive automorphism $\theta$ of the Lie algebra $\g$ and the decomposition
\begin{equation}\label{eq:lie}
  \g = \k \oplus \p, \qquad
  \text{where }\g = \operatorname{Lie}(\G),
  \quad \k = \operatorname{Lie}(\K),
\end{equation}
with $\k$ and $\p$ the $+1$ and $-1$ eigenspaces of $\theta$, respectively \cite{Helgason2001}. In particular, $\K$ is a compact Lie group. The adjoint representation of $\K$ on $\g$ naturally restricts to a representation on $\p$, known as the \emph{isotropy representation} of $\K$. In this case, there exist an $r$-dimensional maximal abelian subspace $\a\subseteq\p$ and a closed convex cone $\a_+\subseteq\a$, known as a closed Weyl chamber, such that every element in $\p$
can be obtained from $\K$'s action on $\a_+$ \cite{ten,Tam1999,Knapp2002}. The resulting factorization is called an
\emph{isotropy decomposition}which is also discussed in the context of Eaton triples \cite{Tam1999} and normal decomposition systems \cite{lewis2000convex}.

\begin{theorem}[Isotropy decomposition]\label{thm:decomp}
For any $A \in \p$, there exists a unique element $A_+ \in \a_+$ such that $A = U\cdot A_+$ for some $U\in\K$. Moreover, for all $Y \in \a_+$,
\begin{equation}\label{eq:var}
  \langle A_+, Y \rangle
  = \max_{U \in \K}\,
  \langle U\cdot A, Y \rangle.
\end{equation}
\end{theorem}

Fix a linear isomorphism $\rho_0:\a\to\R^r$. It induces a map $\rho: \p \to \R^r$ defined by $\rho(A) = \rho_0(A_+)$.
For the purpose of this article, we consider the seven isotropy representations listed in Table~\ref{tab:decomp}. While the first four are the standard eigenvalue and singular value decompositions over $\R$ and $\C$, the last three deserve a few words.

Over the isotropy representation $\Sym(\C^n)$ with congruence action by $\U(n)$, the closed Weyl chamber is $\a_+=\{\diag(v):v_1\geq\cdots\geq v_n\geq0\}$, and for all $A \in \Sym(\C^n)$, there exists $U \in \U(n)$ such that $A = U \diag(\sigma(A)) U^\tp$. This is the Autonne--Takagi decomposition \cite{Autonne15,Takagi1924}.

Over the isotropy representation $\Alt(\R^{2n})$ with congruence action by $\O(2n)$, the closed Weyl chamber is $\a_+=\{\diag(v)\otimes J: v_1\geq\cdots\geq v_n\geq0\}$. Since $A$ is real and skew-symmetric, its eigenvalues occur in pairs $\pm \ii s_i$, where $s_i \geq0$ are singular values of $A$. Equivalently, its singular values occur in equal pairs,
$\sigma_{2i-1}(A)=\sigma_{2i}(A) = s_i$ for $i=1,\dots,n$. Let $\rho(A)= \begin{bmatrix}
    s_1 & \dots & s_n
\end{bmatrix}^\tp$. For all $A \in \Alt(\R^{2n})$, there exists $U \in \O(2n)$ such that $A = U(\diag(\rho(A))\otimes J)U^\tp$. This is the real Youla decomposition \cite{MW1931}, and its complex analogue over $\Alt(\C^{2n})$ has the same form, with $U\in\U(2n)$ \cite{Youla}.

\begin{table*}[htb]
\centering
\tabulinesep=0.75ex
\resizebox{\textwidth}{!}{%
\begin{tabu}{@{}l|l|l|l|l|l|l@{}}
$\p$ & $\K$ & $\a_+$ & $\rho(A)$
& $g\cdot A$ & decomposition & matrix form\\
\hline

$\Herm(\R^n)$
& $\O(n)$
& $\{\diag(v): v_1 \geq \cdots \geq v_n\}$
& $\lambda(A)$
& $UAU^\tp$
& real EVD
& $A = U\Lambda U^\tp$\\

$\Herm(\C^n)$
& $\U(n)$
& $\{\diag(v): v_1 \geq \cdots \geq v_n\}$
& $\lambda(A)$
& $UAU^*$
& complex EVD
& $A = U\Lambda U^*$\\

$\R^{m \times n}$
& $\O(m) \times \O(n)$
& $\{\diag_{m,n}(v): v_1 \geq \cdots \geq v_n \geq 0\}$
& $\sigma(A)$
& $UAV^\tp$
& real SVD
& $A = U\Sigma V^\tp$\\

$\C^{m \times n}$
& $\U(m) \times \U(n)$
& $\{\diag_{m,n}(v): v_1 \geq \cdots \geq v_n \geq 0\}$
& $\sigma(A)$
& $UAV^*$
& complex SVD
& $A = U\Sigma V^*$\\

$\Sym(\C^n)$
& $\U(n)$
& $\{\diag(v): v_1 \geq \cdots \geq v_n \geq 0\}$
& $\sigma(A)$
& $UAU^\tp$
& Autonne--Takagi
& $A = U\Sigma U^\tp$\\

$\Alt(\R^{2n})$
& $\O(2n)$
& $\{\diag(v)\otimes J: v_1 \geq \cdots \geq v_n \geq 0\}$
& $\bigl(\sigma_{2i}(A)\bigr)_{i=1}^n$
& $UAU^\tp$
& real Youla
& $A = U(\Sigma \otimes J)U^\tp$\\

$\Alt(\C^{2n})$
& $\U(2n)$
& $\{\diag(v)\otimes J: v_1 \geq \cdots \geq v_n \geq 0\}$
& $\bigl(\sigma_{2i}(A)\bigr)_{i=1}^n$
& $UAU^\tp$
& complex Youla
& $A = U(\Sigma \otimes J)U^\tp$
\end{tabu}%
}
\caption{Examples of isotropy decompositions.}
\label{tab:decomp}
\end{table*}

The eigenvalue, singular value, and Autonne--Takagi decompositions also admit quaternionic analogues \cite{HZ1012}, which we omit for simplicity of presentation. Together with the seven decompositions in Table~\ref{tab:decomp}, these give the ten classical families of isotropy decompositions \cite{ten}. By the classification of classical symmetric spaces,
all isotropy decompositions arising from these spaces can be obtained from the ten families and their direct products, up to trivial factors. In the rest of the article, we reserve the notation $\p$ for the seven vector spaces listed in Table~\ref{tab:decomp}.

\section{Reduction via isotropic divergences}
\subsection{Isotropic divergence}
Recall that for a manifold $\mathcal{M}$, a \emph{divergence} on $\Omega \subseteq \mathcal{M}$ is a $C^2$ function $D:\Omega\times\Omega\to[0,\infty)$ such that 
$D(x,y)=0$ if and only if $x=y$, and in local coordinates,
\begin{equation}\label{eq:taylor}
    D(x,x+h)=\frac12\, g_x(h,h)+o(\lVert h\rVert^2),
\end{equation}
where $g_x$ is a positive-definite quadratic form on $T_x\mathcal{M}$ \cite[Definition~1.1]{Amari2016}. In this article, we focus on matrix divergences with with $\mathcal{M}=\p$ and $\Omega\subseteq\p$ open and $\K$-invariant, i.e., $\K \cdot \Omega \subseteq \Omega$. 
\begin{definition}
A divergence $D$ on $\Omega \subseteq \p$ is \emph{$\K$-invariant} if $D(X,Y) = D(U \cdot X, U \cdot Y)$ for all $U \in \K$ and $X$, $Y \in \Omega$. We say $D$ is $\K$-\emph{isotropic} if it is $\K$-invariant and
\[
\min_{U \in \K} D(U \cdot A, Y) = D(A_+, Y)
\]
for all $A \in \Omega$ and $Y \in \Omega \cap \a_+$.    
\end{definition}

Equivalently, $\K$-isotropy can be viewed as a divergence analogue of the $1$-Lipschitz property for the map $A\mapsto A_+$, generalizing \eqref{eq:ing2}.

\begin{proposition}\label{prop:char}
Let $\Omega\subseteq\p$ be $\K$-invariant and open and $D$ be a $\K$-invariant divergence on $\Omega$. Then the following are equivalent: 
\begin{enumerate}[\normalfont(i)]
    \item For all $A,B\in\Omega$, \label{item:char1}
    \[
    D(A_+,B_+)\leq D(A,B).
    \]
    \item For all $A\in\Omega$ and $Y\in\Omega\cap\a_+$, \label{item:char2}
    \[
        \min_{U\in K}D(U \cdot A,Y)=D(A_+,Y).
    \]
    \item For all $A\in\Omega$ and $Y\in\Omega\cap\a_+$, \label{item:char3}
    \[
    \min_{U\in K}D(Y,U \cdot A)=D(Y,A_+).
    \]
\end{enumerate}
\end{proposition}

\begin{proof}
We first prove that \ref{item:char1} implies \ref{item:char2}. Let $A\in\Omega$ and $Y\in\Omega\cap\a_+$. By Theorem~\ref{thm:decomp} and the uniqueness of $A_+$ and $Y_+$, $(U \cdot A)_+=A_+$ for all $U \in \K$ and $Y_+=Y$. Using the assumption \ref{item:char1}, we obtain that
\[
D(A_+,Y)
=
D\bigl((U \cdot A)_+,Y_+\bigr)
\leq
D(U \cdot A,Y)
\]
for all $U\in \K$. On the other hand, let $U_A\in K$ such that $U_A \cdot A=A_+$, so $D(U_A \cdot A,Y)=D(A_+,Y)$. Since the lower bound is attained, we obtain \ref{item:char2}. The implication from \ref{item:char1} to \ref{item:char3} follows analogously.

Now we show that \ref{item:char2} implies \ref{item:char1}. Let $A,B\in\Omega$. Choose $U_B \in K$ such that $B= U_B \cdot B_+$ by Theorem~\ref{thm:decomp}. Since $B_+\in\a_+$, \ref{item:char2}, by \ref{item:char2} and the $\K$-invariance of $D$, we obtain
\[
D(A,B) = 
D\bigl(U_B^{-1} \cdot A,B_+\bigr)
\geq
D\bigl((U_B^{-1} \cdot A)_+,B_+\bigr) = D(A_+,B_+),
\]
where the last equality follows from \eqref{eq:var} and the uniqueness of $A_+$, which gives \ref{item:char1}. The implication from \ref{item:char3} to \ref{item:char1} follows from a verbatim argument.
\end{proof}

\subsection{Reduced divergence}

We now formally state and prove Theorem~\ref{thm:intro-reduction}. 

\begin{theorem}[Reduced divergence]\label{thm:reduce}
Let $\mathcal{S}\subseteq\Omega\subseteq\p$ such that $\mathcal{S}$ is $\K$-invariant and $\Omega$ is open and $\K$-invariant. If $D$ is a $\K$-isotropic divergence on $\Omega$, then there exists a unique function $\widehat D$ on $\rho(\Omega)$, which is the
restriction of a divergence on an open neighborhood of $\rho(\Omega)$, such that
\begin{equation}\label{eq:equiv}
    D(A_+,B_+)
=
\widehat D\bigl(\rho(A),\rho(B)\bigr)
\qquad
\text{for all }A,B\in\Omega.
\end{equation}
Moreover, for every $A\in\Omega$,
\begin{equation}\label{eq:right}
    \inf_{X\in\mathcal{S}}D(A,X)
=
\inf_{Y\in\mathcal{S}\cap\a_+}D(A_+,Y)
=
\inf_{v\in\rho(\mathcal{S})}
\widehat D\bigl(\rho(A),v\bigr),
\end{equation}
and
\begin{equation}\label{eq:left}
    \inf_{X\in\mathcal{S}}D(X,A)
=
\inf_{Y\in\mathcal{S}\cap\a_+}D(Y,A_+)
=
\inf_{v\in\rho(\mathcal{S})}
\widehat D\bigl(v,\rho(A)\bigr).
\end{equation}
\end{theorem}
\begin{proof}
By \cite[ChapterVII, Proposition~2.2 and Theorem~2.22]{Helgason2001}, $(\K \cdot A) \cap \a_+$ is a singleton for all $A \in \p$. Since $\Omega$ is $\K$-invariant, $A_+\in\Omega\cap\a_+$ whenever $A \in \Omega$. We hence obtain that $\rho(\Omega)=\rho(\Omega\cap\a_+)$. Moreover, for $X$, $Y\in\Omega\cap\a_+$, since $X_+=X$, $Y_+=Y$, and $\rho_0:\a\to\R^r$ is an isomorphism,
\begin{equation}\label{eq:reduce1}
    \rho(X)=\rho(Y)
\iff
\rho_0(X)=\rho_0(Y)
\iff
X=Y.
\end{equation}
Consequently, the map $\rho \vert_{\Omega\cap\a_+}: \Omega\cap\a_+\to\rho(\Omega)$ is a bijection. Consider
\[
\widehat{D}: \rho(\Omega)\times\rho(\Omega)\to[0,\infty),
\qquad
\widehat{D}(x,y)
\coloneqq
D\bigl(
\rho\vert_{\Omega\cap\a_+}^{-1}(x),
\rho\vert_{\Omega\cap\a_+}^{-1}(y)
\bigr).
\]
Clearly, $\widehat{D}$ is well defined by the bijectivity of $\rho\vert_{\Omega\cap\a_+}$. We first show that $\widehat D$ is a divergence. Since $\Omega$ is open in $\p$ and $\rho_0:\a\to\R^r$ is a linear isomorphism, the image $\rho_0(\Omega\cap\a)$ is open in $\R^r$. Define
\[
D_{\a}:\rho_0(\Omega\cap\a) \times \rho_0(\Omega\cap\a) \to [0,\infty),
\qquad
D_{\a}(x,y)
\coloneqq
D\bigl(\rho_0^{-1}(x),\rho_0^{-1}(y)\bigr).
\]
Since $D$ is $C^2$ and $\rho_0^{-1}$ is linear, the composition $D_{\a}$ is $C^2$. Moreover, since $\rho(\Omega)\subseteq\rho_0(\Omega\cap\a)$ and $\rho(X)=\rho_0(X)$ for $X\in\a_+$, we have that $\widehat D =D_{\a}\big \vert_{\rho(\Omega)\times\rho(\Omega)}$. Consequently, $\widehat D$ is the restriction of a $C^2$ function defined on an open neighborhood of $\rho(\Omega)$.

The nonnegativity of $\widehat{D}$ follows trivially from the nonnegativity of $D$. For $x,y\in\rho(\Omega)$, we have
\[
\widehat D(x,y)=0
\iff
D\bigl(\rho\vert_{\Omega\cap\a_+}^{-1}(x),\rho\vert_{\Omega\cap\a_+}^{-1}(y)\bigr)=0
\iff
\rho\vert_{\Omega\cap\a_+}^{-1}(x)=\rho\vert_{\Omega\cap\a_+}^{-1}(y)
\iff
x=y.
\]

It remains to verify \eqref{eq:taylor}. Let $x\in\rho(\Omega)$. Since $D$ is a divergence, there exists a positive-definite quadratic form $g_{\rho\vert_{\Omega\cap\a_+}^{-1}(x)}$ on $T_{\rho\vert_{\Omega\cap\a_+}^{-1}(x)}\p$ such that
\[
D(\rho\vert_{\Omega\cap\a_+}^{-1}(x),\rho\vert_{\Omega\cap\a_+}^{-1}(x)+H)
=
\frac12 g_{\rho\vert_{\Omega\cap\a_+}^{-1}(x)}(H,H)+o(\lVert H\rVert^2).
\]
Let $X=\rho_0^{-1}(x)$. For any sufficiently small $h\in\R^r$ such that $x+h\in\rho(\Omega)$, we have $\rho\vert_{\Omega\cap\a_+}^{-1}(x+h)
=
X+\rho_0^{-1}(h)$, and hence
\[
\widehat D(x,x+h)=
D\bigl(X,X+\rho_0^{-1}(h)\bigr)=
\frac12
g_X\bigl(\rho_0^{-1}(h),\rho_0^{-1}(h)\bigr)
+o(\lVert h\rVert^2).
\]
We thus obtain that the quadratic form induced by $\widehat D$ at $x$ is
\[
\widehat g_x(h,k)
\coloneqq
g_{\rho\vert_{\Omega\cap\a_+}^{-1}(x)}\bigl(\rho_0^{-1}(h),\rho_0^{-1}(k)\bigr),
\]
which is positive definite because $g_{\rho\vert_{\Omega\cap\a_+}^{-1}(x)}$ is positive definite and $\rho_0^{-1}$ is injective. We have thus shown that $\widehat D$ is a divergence on $\rho(\Omega)$.

We now show uniqueness. Let $A,B\in\Omega$. By construction, we have $\rho\vert_{\Omega\cap\a_+}^{-1}(\rho(A))=A_+$ and $\rho\vert_{\Omega\cap\a_+}^{-1}(\rho(B))=B_+$ and hence $\widehat D\bigl(\rho(A),\rho(B)\bigr) = D(A_+,B_+)$. If $\widetilde D$ is another divergence on $\rho(\Omega)$ satisfying $D(A_+,B_+) = \widetilde D\bigl(\rho(A),\rho(B)\bigr)$ for all $A, B \in \Omega$, then for all $x$, $y \in \rho(\Omega)$
\[
\widetilde D(x,y)
=
D\bigl(\rho_0^{-1}(x),\rho_0^{-1}(y)\bigr)
=
\widehat D(x,y).
\]

Lastly, it suffices to prove \eqref{eq:right} as \eqref{eq:left} follows from an analogous argument. By Proposition~\ref{prop:char}, $D(A_+,X_+)\leq D(A,X)$ for all $A, X \in \Omega$. Since $\mathcal S$ is $\K$-invariant, we have that $X_+\in\mathcal S\cap\a_+$ if $X\in\mathcal S$, and as a result,
\[
\inf_{X\in\mathcal S}D(A,X)
\geq
\inf_{Y\in\mathcal S\cap\a_+}D(A_+,Y).
\]
Conversely, we choose $U_A\in\K$ such that $A=U_A \cdot A_+$ using Theorem~\ref{thm:decomp}. For every $Y\in\mathcal S\cap\a_+$, we have $U_A \cdot Y\in\mathcal S$, and $\K$-invariance of $D$ yields that
\[
D\bigl(A,U_A \cdot Y\bigr)
=
D\bigl(U_A \cdot A_+,U_A \cdot Y\bigr)
=
D(A_+,Y).
\]
We hence obtain the first equality in \eqref{eq:right}. 

By \eqref{eq:reduce1}, the map $\rho\vert_{\mathcal S\cap\a_+}: \mathcal S\cap\a_+\to\rho(\mathcal S)$ is a bijection. Applying \eqref{eq:equiv} yields
\[
\inf_{Y\in\mathcal S\cap\a_+}D(A_+,Y)
=
\inf_{v\in\rho(\mathcal S)}
\widehat D\bigl(\rho(A),v\bigr). \qedhere
\]
\end{proof}

As illustrated in the introduction, Theorem~\ref{thm:reduce} allows us to reduce the matrix nearness problem \eqref{eq:near} to a vector nearness problem of approximating $\rho(A)$ over $\widehat{D}$, which we call the \emph{reduced divergence} of $D$. Using isotropic decomposition, we can immediately reproduce a minimizer of the original problem from a minimizer of the reduced problem.

\begin{corollary}\label{cor:recover}
Under the assumptions of Theorem~\ref{thm:reduce}, let $A \in \Omega$ with isotropic decomposition $A = U \cdot A_+$. If 
\[
u_* \in \argmin_{u\in\rho(\mathcal S)}
\widehat D\bigl(\rho(A),u\bigr) \quad \text{and} \quad v_* \in \argmin_{v\in\rho(\mathcal S)}
\widehat D\bigl(v,\rho(A)\bigr)
\]
then 
\[
   X_* \coloneqq U\cdot \rho_0^{-1}(u_*) \in \argmin_{X \in \mathcal{S}} D(A,X) \quad \text{and} \quad Y_* \coloneqq U\cdot \rho_0^{-1}(v_*) \in \argmin_{Y \in \mathcal{S}} D(Y,A).
\]
\end{corollary}
\begin{proof}
Since $u_*$, $v_*\in\rho(\mathcal S)$ and $\mathcal S$ is $\K$-invariant, we have $\rho_0^{-1}(u_*),\rho_0^{-1}(v_*)\in\mathcal S\cap\a_+$. Therefore, $X_*, Y_* \in \mathcal S$. Using the $\K$-invariance of $D$, \eqref{eq:equiv}, \eqref{eq:right}, and \eqref{eq:left}, we obtain
\[
D\bigl(A,X_*\bigr)
=\widehat D\bigl(\rho(A),u_*\bigr)
=\min_{X\in\mathcal S}D(A,X), \quad D\bigl(Y_*,A\bigr)
=\widehat D\bigl(v_*,\rho(A)\bigr)
=\min_{Y\in\mathcal S}D(Y,A).  \qedhere
\]
\end{proof}

\section{Examples of isotropic divergences}\label{sec:divergences}

Isotropic divergences are ubiquitous in various applications. To illustrate this, we now show that all divergences in Table~\ref{tab:isotropic-divergences} are isotropic and derive the corresponding reduced divergences. For convenience, we fix the constants 
\[
c_\p = \begin{cases}
    2 & \p = \Alt(\F^{2n})\\
    1 & \text{otherwise}.
\end{cases}
\]
In practice, there are different classes of divergences such as the Bregman divergence and the Jensen divergence, which are characterized by a given scalar function. In these cases, we provide general sufficient conditions for isotropy.
    
\subsection{Bregman divergences}

Let $\Omega\subseteq\p$ be an open convex set. For a differentiable strictly convex function $\varphi:\Omega\to\R$, the Bregman divergence
$D_\varphi$ is defined by
\begin{equation}\label{eq:bregman}
    D_\varphi(p, q) \coloneqq \varphi(p) - \varphi(q) - \langle \nabla \varphi(q), p - q \rangle.
\end{equation}
The isotropy of Bregman divergence can be characterized by $\varphi$.

\begin{proposition}[Isotropic Bregman divergence]
\label{prop:bregman}
Let $\Omega\subseteq\p$ be open, convex, and $\K$-invariant. If $\varphi \in C^3(\Omega)$ is $\K$-invariant and $\nabla^2\varphi(X)$ is positive definite
for every $X\in\Omega$,
then $D_\varphi$ is $\K$-isotropic.
\end{proposition}

\begin{proof}
Since $\nabla^2\varphi$ is positive definite, $\varphi$ is strictly convex. So $D_\varphi$is nonnegative and $D_\varphi(X,Y) = 0$ if and only if $X=Y$. Moreover, since $\varphi\in C^3(\Omega)$, $D_\varphi\in C^2(\Omega\times\Omega)$, and
\[
D_\varphi(X,X+H)
=
\frac12\langle\nabla^2\varphi(X)H,H\rangle
+o(\lVert H \rVert^2).
\]
Thus $D_\varphi$ is a divergence. Its $\K$-invariance follows from that of $\varphi$.

Let $A\in\Omega$ and $Y\in\Omega\cap\a_+$. Since $\varphi$ is $\K$-invariant, for all $U \in \K$,
\begin{equation}\label{eq:berg-decomp}
D_\varphi(U \cdot A,Y)
=
\varphi(A)-\varphi(Y)
+\langle\nabla\varphi(Y),Y\rangle
-\langle\nabla\varphi(Y),U \cdot A\rangle.
\end{equation}

Since the first three terms do not depend on $U$, minimizing $D_\varphi(U \cdot A,Y)$ is equivalent to maximizing $\langle\nabla\varphi(Y),U \cdot A\rangle$. Since $\varphi$ is strictly convex, applying \cite[Corollary~4.7]{Lewis1996} yields
\[
\nabla\varphi(Y)
=
\nabla\varphi(Y_+)
=
\bigl(\nabla\varphi(Y)\bigr)_+
\in\a_+
\]
for every $Y\in\Omega\cap\a_+$. By Theorem~\ref{thm:decomp},
\[
\max_{U\in\K}
\langle\nabla\varphi(Y), U \cdot A\rangle
=
\langle\bigl(\nabla\varphi(Y)\bigr)_+,A_+\rangle
=
\langle\nabla\varphi(Y),A_+\rangle.
\]
Substituting it into \eqref{eq:berg-decomp} yields
\begin{equation}\label{eq:prop41}
    \min_{U\in\K}D_\varphi(U \cdot A,Y)
=
\varphi(A)-\varphi(Y)
+\langle\nabla\varphi(Y),Y\rangle
-\langle\nabla\varphi(Y),A_+\rangle
=
D_\varphi(A_+,Y). \qedhere
\end{equation}
\end{proof}

Proposition~\ref{prop:bregman} allows us to verify the $\K$-isotropy of a wide range of familiar matrix Bregman divergences, starting with the simplest example of the \emph{Frobenius divergence}.

\subsubsection{Frobenius divergence}
Let $\Omega=\p$, where $\p$ is any isotropy representation in Table~\ref{tab:decomp}. The squared Frobenius norm $\varphi(X)=\frac12\lVert X\rVert_{\mathsf F}^2$ induces the Bregman divergence 
\[
D_{\mathsf F}(X,Y)=\frac12\lVert X-Y\rVert_{\mathsf F}^2 \qquad \text{for all } X, Y \in \p.
\]
Clearly, $\varphi$ is smooth, strictly convex, and $\K$-invariant. Applying Theorem~\ref{thm:reduce} thus yields the reduced divergence 
\[
\begin{aligned}
\widehat D_{\mathsf F}(u,v) =D_{\mathsf F}\bigl(
\rho\vert_{\a_+}^{-1}(u),
\rho\vert_{\a_+}^{-1}(v)\bigr) =\frac12\lVert\rho_0^{-1}(u-v)\rVert_{\mathsf F}^2 
=\frac{c_{\p}}2\lVert u-v\rVert_2^2 \qquad \text{for all } u, v \in \rho(\p).
\end{aligned}
\]

In the rest of this section, the reduced divergences are all similarly derived using Theorem~\ref{thm:reduce} with straightforward calculations once we verify the $\K$-isotropy of $D$. For clarity of the exposition, we will hence omit the intermediate calculations and directly present the formula of the reduced divergences. 

\subsubsection{Von Neumann divergence}\label{sec:vN}
Another commonly used divergence is the \emph{von Neumann divergence}, also known as quantum relative entropy, which arise in the literature of quantum mechanics \cite{DT2007,NielsonChuang2000}. Let $\p=\Herm(\mathbb F^n)$ and $\Omega=\Herm_+(\mathbb F^n)$. For $X \in \Herm_+(\mathbb F^n)$, the von Neumann entropy $\varphi(X)=\tr(X\log X-X)$ induces it as the Bregman divergence
\[
D_{\mathsf{vN}}(X,Y)
=\tr\bigl(X(\log X-\log Y)-X+Y\bigr) \qquad \text{for all } X, Y \in \Herm_+(\mathbb F^n).
\]
The von Neumann entropy $\varphi$ is clearly smooth and $\K$-invariant. For every $0 \neq H \in \Herm(\F^n)$, a direct computation shows that
\[
\langle H,\nabla^2\varphi(X)[H]\rangle
=
\int_0^\infty
\bigl\|(X+tI)^{-1/2}H(X+tI)^{-1/2}\bigr\|_{\mathsf F}^2
\,dt>0.
\]
Therefore, $\nabla^2\varphi(X)$ is positive definite and hence $D_{\mathsf{vN}}$ is $\K$-isotropic by Proposition~\ref{prop:bregman}. Its reduced divergence is given by

\[
\widehat D_{\mathsf{vN}}(u,v)
=\sum_{i=1}^n
\biggl(u_i\log\frac{u_i}{v_i}-u_i+v_i\biggr) \qquad \text{for all } u, v \in \rho(\Omega).
\]

\subsubsection{Logdet divergence}
The domain of positive definite matrices also arise naturally as covariance matrices in statistics and learning theory. In this context, a natural divergence is the \emph{logdet divergence}, also known as Stein's loss in the statistics literature \cite{logdet1,logdet2}. Let $\p=\Herm(\mathbb F^n)$ and $\Omega=\Herm_+(\mathbb F^n)$. The logdet function $\varphi(X)=-\log\det X$ induces the logdet divergence 
\[
D_{\mathsf{LD}}(X,Y)
=\tr(XY^{-1})-\log\det(XY^{-1})-n \qquad \text{for all } X, Y \in \Herm_+(\mathbb F^n).
\]
Notice that for $X,Y\in\Herm_+(\mathbb R^n)$,
\[
D_{\mathsf{LD}}(X,Y) = 2 D_{\mathrm{KL}}\bigl(
\mathcal N(0,X)\,\Vert \,\mathcal N(0,Y)
\bigr),
\]
where $D_{\mathrm{KL}}$ is the Kullback--Leibler divergence \cite{KL}. Since the Hessian of $\varphi$ is
\begin{equation}\label{eq:logdet-hess}
    \nabla^2\varphi(X)[H] = X^{-1}HX^{-1} \qquad \text{for all } X \in \Herm_+(\F^n)
\end{equation}
and $\K$-invariant by $\K$-invariance of determinant, the logdet divergence is isotropic by Proposition~\ref{prop:bregman} and the corresponding reduced divergence is
\[
\widehat D_{\mathsf{LD}}(u,v)
=\sum_{i=1}^n
\biggl(\frac{u_i}{v_i}-\log\frac{u_i}{v_i}-1\biggr) \qquad \text{for all } u, v \in \rho(\Omega).
\]

\subsubsection{Bregman divergence in canonical correlation analysis}
Another closely related domain is
\[
\p_{<\mu} \coloneqq \{X \in \p: \mu I-X^*X \succ 0\}.
\]
In particular, $\p_{< 1}$ is the domain of strictly contractive matrices, which relates to the positive definite matrices as
\begin{equation}\label{eq:pd}
X \in \p_{<1} \quad \Longleftrightarrow \quad \mathcal{C}(X) \coloneqq \begin{bmatrix}
    I &X^* \\ X &I
\end{bmatrix} \succ 0.
\end{equation}
From the statistical point-of-view, they arise naturall in the context of canonical correlation analysis, giving the following divergence \cite{Hotelling1936,Anderson2003,BachJordan2005}. Let $\p=\mathbb R^{m\times n}$ and $\Omega=\R^{m\times n}_{<1}$. For $X, Y \in \R^{m\times n}_{<1}$, consider  
\[
\begin{aligned}
D_{\mathsf{GKL}}(X,Y)
&\coloneqq \frac12\log\frac{\det(I-Y^*Y)}{\det(I-X^*X)} +\Re\tr\bigl[(I-Y^*Y)^{-1}Y^*(Y-X)\bigr]\\
&=D_{\mathrm{KL}}(N(0,\mathcal C(X))\Vert N(0,\mathcal C(Y))).
\end{aligned}
\]
It is the Bregman divergence induced by
\[
\varphi(X)
=
-\frac12\log\det(I-X^\tp X)
=
-\frac12\log\det\mathcal C(X).
\]
Clearly, $\varphi$ is smooth and invariant under the action of $\O(m)\times\O(n)$. Since $\mathcal C(X)$ is positive definite affine in $X$, $\nabla^2\varphi(X)$ is positive definite for all $X \in \Omega$ by \eqref{eq:logdet-hess}. By Proposition~\ref{prop:bregman}, $D_{\mathsf{GKL}}$ is $\K$-isotropic. By Theorem~\ref{thm:reduce}, its reduced divergence is
\[
\widehat D_{\mathsf{GKL}}(u,v)
=\frac12\sum_{i=1}^r\Biggl[
\log\frac{1-v_i^2}{1-u_i^2}
-\frac{2v_i(u_i-v_i)}{1-v_i^2}\Biggr].
\]

\subsubsection{Spectral hypentropy}
Lastly, another way of constructing isotropic Bregman divergence is to apply a strictly convex function $f$ component-wise to $\rho_i(X)$'s. The spectral hypentropy, introduced in the context of machine learning, is such an example \cite{hypentropy}. Let $\Omega = \p = \Herm(\F^n)$. For $\beta>0$,
let
\[
D_\beta(X,Y)\coloneqq\tr\biggl[X\biggl(\operatorname{arsinh}\frac{X}{\beta}-\operatorname{arsinh}\frac{Y}{\beta}\biggr)+\sqrt{Y^2+\beta^2 I}-\sqrt{X^2+\beta^2 I} \biggr] \qquad \text{for all } X, Y \in \p.
\]
It is the Bregman divergence induced by $\varphi_\beta(X)=\sum_{i=1}^r f_\beta(\lambda_i(X))$ where $f_\beta(t) =t\operatorname{arsinh}(t/\beta) -\sqrt{t^2+\beta^2}+\beta$. Since $f_\beta$ is smooth and
$f_\beta''(t)=(t^2+\beta^2)^{-1/2}>0$, a direct computation shows that $\varphi_\beta$ is smooth with positive-definite Hessian. Combined with the clear fact that $\varphi_\beta$ is $\K$-invariant, Proposition~\ref{prop:bregman} yields that $D_\beta$ is $\K$-isotropic with reduced divergence
\[
\widehat D_\beta(u,v)=\sum_{i=1}^r\biggl[u_i\biggl(
\operatorname{arsinh}\frac{u_i}{\beta}
-\operatorname{arsinh}\frac{v_i}{\beta}\biggr)-\sqrt{u_i^2+\beta^2}+\sqrt{v_i^2+\beta^2}
\biggr] \qquad \text{for all } u,v \in \rho(\p).
\]

\subsection{Jensen divergences}
For an open convex set $\Omega\subseteq\p$ and a strictly convex function $\varphi:\Omega\to\mathbb R$, the Jensen divergence, also known as Burbea--Rao divergences \cite{NielsenBoltz2011}, is defined by   
\begin{equation}\label{eq:jensen}
    J_\varphi(X,Y)
=\frac{\varphi(X)+\varphi(Y)}2
-\varphi\left(\frac{X+Y}{2}\right) \qquad \text{for all } X, Y \in \Omega.
\end{equation}
Notice that it is always symmetric. 

\begin{proposition}[Isotropy of Jensen divergences]
\label{prop:jensen-isotropy}
Let $\Omega\subseteq\p$ be open, convex, and
$\K$-invariant. If $\varphi\in C^2(\Omega)$ is $\K$-invariant and $\nabla^2\varphi(X)$ is positive definite for every $X\in\Omega$,
then $J_\varphi$ is $\K$-isotropic.
\end{proposition}

\begin{proof}
Let $X,Y\in \Omega$. By \eqref{eq:jensen}, $J_\varphi$ is $\K$-invariant. Since $\nabla^2\varphi$ is positive definite,
$\varphi$ is strictly convex, so $J_\varphi$ is nonnegative and $J_\varphi(X,Y) = 0$ if and only if $X = Y$. Moreover, $J_\varphi\in C^2(\Omega\times\Omega)$, and
\[
J_\varphi(X,X+H)
=
\frac18\langle\nabla^2\varphi(X)H,H\rangle
+o(\|H\|^2).
\]
Thus $J_\varphi$ is a divergence. For convenience of notation, $M = (X+Y)/2$. Then by \eqref{eq:jensen} and \eqref{eq:bregman},
\[
J_\varphi(X,Y)
=
\frac12\bigl(
D_\varphi(X,M)+D_\varphi(Y,M)
\bigr).
\]
By \eqref{eq:prop41}, $D_\varphi(A_+,B_+)\leq D_\varphi(A,B)$ for all $A,B\in\Omega$. Consequently,
\begin{align*}
J_\varphi(X,Y)
&\geq
\frac12\bigl(
D_\varphi(X_+,M_+)
+
D_\varphi(Y_+,M_+)
\bigr)\\
&=
J_\varphi(X_+,Y_+)
+
D_\varphi\biggl(
\frac{X_++Y_+}{2},M_+
\biggr) \geq J_\varphi(X_+,Y_+). \qedhere
\end{align*}
\end{proof}

\subsubsection{$S$-divergence}
The S-divergence, also known as the Jensen--Bregman logdet divergence, arises in a variety of applications including machine learning and statistics \cite{Cherian2013,Sra2016}. Let $\p = \Herm(\F^n)$ and $\Omega = \Herm_+(\F^n)$, the S-divergence is the Jensen divergence induced by $\varphi(X) = -\log \det (X)$:
\[
D_\mathsf{S}(X,Y) = \log \det \biggl( \frac{X+Y}{2} \biggr)
 - \frac{1}{2}\log \det (XY) \qquad \text{for all } X, Y \in \Herm_+(\mathbb F^n).
\]
Since $\varphi$ is smooth, $\K$-invariant, and has positive-definite Hessian by
\eqref{eq:logdet-hess}, Proposition~\ref{prop:jensen-isotropy} and Theorem~\ref{thm:reduce} yields its reduced divergence
\[
\widehat D_{\mathsf S}(u,v)
=\sum_{i=1}^n
\log\frac{u_i+v_i}{2\sqrt{u_iv_i}} \qquad \text{for all } u,v \in \rho(\Omega).
\]

\subsubsection{Quantum Jensen--Shannon divergence}
Another exquisite example of Jensen divergence is the quantum Jensen--Shannon divergence \cite{Majtey2005,Sra2021}, induced by the von Neumann entropy on  $\p=\Herm(\mathbb F^n)$ and $\Omega=\Herm_+(\mathbb F^n)$:
\[
D_{\mathsf{QJS}}(X,Y)
=\frac12\tr(X\log X+Y\log Y)
-\tr\biggl(\frac{X+Y}{2}\log\frac{X+Y}{2}\biggr) \qquad \text{for all } X, Y \in \Herm_+(\mathbb F^n).
\]
As established in Section~\ref{sec:vN}, $\varphi$ is smooth, $\K$-invariant, and has positive-definite Hessian. We apply Proposition~\ref{prop:jensen-isotropy} and Theorem~\ref{thm:reduce} to obtain its reduced divergence
\[
\widehat D_{\mathsf{QJS}}(u,v)
=\sum_{i=1}^n\left[
\frac{u_i\log u_i+v_i\log v_i}{2}
-\frac{u_i+v_i}{2}\log\frac{u_i+v_i}{2}
\right] \qquad \text{for all } u, v \in \rho(\Omega).
\]

\subsection{Calabi's diastases}

Let $(\mathcal M,g)$ be a real-analytic K\"ahler manifold with local K\"ahler potential $\Phi$. Calabi's diastasis \cite{Calabi1953} is the locally defined function
\[
C(p,q)
=\Phi(p,\overline p)+\Phi(q,\overline q)
-\Phi(p,\overline q)-\Phi(q,\overline p),
\]
where the mixed arguments are defined via analytic continuation. In the case of
$\mathcal{M} = \p_{<\mu}$ for $\F = \C$, it provides meaningful divergences with meaningful applications.

\begin{proposition}[Diastasis functions as isotropic divergences]
\label{prop:calabi}
Let $\p=\C^{m\times n}$, $\Sym(\C^n)$, or $\Alt(\C^{2n})$ and $\Omega \subseteq \p_{<\mu}$ be open and $\K$-invariant. Fix $\mu\in(0,\infty]$ and let
\[
    f(z)=\sum_{k=1}^{\infty}a_kz^k,
    \qquad a_1>0,\quad a_k\ge0
    \text{ for } k = 2,3,\dots
\]
be holomorphic on $\lvert z \rvert <\mu$. Then $\tr f(X^*X)$ defines a real-analytic K\"ahler potential on $\Omega$, whose Calabi diastasis
\begin{equation}\label{eq:spectral-calabi}
    C_f(X,Y)
    =\tr f(X^*X)+\tr f(Y^*Y)
    -2\Re[\tr f(Y^*X)]
\end{equation}
is a $\K$-isotropic divergence with reduced divergence
\begin{equation}\label{eq:spectral-calabi-reduced}
    \widehat D_f(u,v)
    =c_{\p}\sum_{i=1}^{r}\sum_{k=1}^{\infty}
        a_k(u_i^k-v_i^k)^2,
    \quad \text{for all } u,v \in \rho(\Omega).
\end{equation}
\end{proposition}

\begin{proof}
Let $\Phi_f(X)=\tr f(X^*X)$. Since $\lVert Y^*X\rVert_{\mathsf{op}} \leq\lVert Y\rVert_{\mathsf{op}} \lVert X\rVert_{\mathsf{op}}<\mu$ for all $X,Y\in\Omega$, the power series
\[
\phi_f(X,\overline{Y})
=\sum_{k=1}^\infty a_k
\tr((Y^*X)^k)
=\tr f(Y^*X)
\]
converges.

For every $k\geq1$, by H\"older's inequality for Schatten norms and
\cite[Theorem~3.3.4]{HJ1991},
\begin{equation}\label{eq:calabi1}
\Re\bigl[\tr((Y^*X)^k)\bigr]
\leq \sum_{i=1}^{c_{\p}r}\sigma_i(Y^*X)^k
\leq \sum_{i=1}^{c_{\p}r}
\sigma_i(X)^k\sigma_i(Y)^k
=\tr((Y_+^*X_+)^k).
\end{equation}    
Consequently,
\[
\begin{aligned}
C_f(X,Y)
&=\sum_{k=1}^{\infty}a_k
\bigl[
  \tr((X^*X)^k)
  +\tr((Y^*Y)^k)
  -2\Re\tr((Y^*X)^k)
\bigr]\\
&\geq\sum_{k=1}^{\infty}a_k
\sum_j(\sigma_j(X)^k-\sigma_j(Y)^k)^2
\geq0.
\end{aligned}
\]
In particular,
\begin{equation}\label{eq:bound1}
    C_f(X,Y)\geq a_1\lVert X-Y\rVert_F^2.
\end{equation}

It remains to verify that $\Phi_f$ defines a K\"ahler metric. For $X\in\Omega$ and $H\in\p$, Taylor expansion gives
\[
C_f(X,X+tH)
=|t|^2
\frac{\partial^2}{\partial s\,\partial\overline{s}}
\Phi_f(X+sH)\Big|_{s=0}
+o(|t|^2).
\]
By \eqref{eq:bound1}, dividing by $\lvert t\rvert^2$ and letting $t\to0$ yields
\[
\frac{\partial^2}{\partial s\,\partial\overline{s}}
\Phi_f(X+sH)\Big|_{s=0}
\geq a_1\lVert H\rVert_{\mathsf F}^2.
\]
Therefore, $\Phi_f$ defines a real-analytic K\"ahler metric on $\Omega$ and $C_f$
is its Calabi diastasis. Since $C_f(X,Y)=0$ if and only if $X=Y$ by \eqref{eq:bound1}, $C_f$ is a divergence.

By the unitary invariance of trace, $C_f$ is $\K$-invariant. By \eqref{eq:calabi1}, $C_f(X,Y)\geq C_f(X_+,Y_+)$ for all $X,Y\in\Omega$. Therefore, $C_f$ is $\K$-isotropic and a direct computation yields \eqref{eq:spectral-calabi-reduced}.
\end{proof}

\subsubsection{Calabi's diastasis in bosonic squeezing}

The first example of isotropic diastasis divergence is motivated by bosonic squeezing \cite{BianchiHacklYokomizo2015} in the case of $\p = \Sym(\C^n)$, where a matrix $X \in \Sym(\C^n)_{<1}$ parametrizes a normalized pure bosonic squeezed-vacuum state by
\[
    \ket{\psi_X}
    =\det(I-X^*X)^{1/4}
    \exp\Biggl(\frac12\sum_{j,k=1}^n
            X_{jk}a_j^\dagger a_k^\dagger \Biggr) \ket{0}.
\]
For $X, Y \in \Sym(\C^n)_{<1} = \Omega$, the fidelity between the corresponding states is
\[
    |\langle\psi_X|\psi_Y\rangle|^2
    =\frac{\sqrt{\det(I-X^*X)\det(I-Y^*Y)}}
            {|\det(I-X^*Y)|}.
\]
Consider the divergence defined by
\[
    C_{\mathsf{B}}(X,Y)
    \coloneqq -\log\lvert \braket{\psi_X|\psi_Y} \rvert^2 \qquad \text{for all } X, Y \in \Sym_{<1}(\C^n).
\]
Clearly, it is the Calabi's diastasis induced by 
\[
    f(z)=-\frac12\log(1-z) = \sum_{k=1}^{\infty}\frac{z^k}{2k}, \quad \Phi_f(X)=-\frac12\log\det(I-X^*X).
\]
Using Proposition~\ref{prop:calabi}, we obtain that $C_{\mathsf{B}}$ is a $\K$-isotropic divergence with reduced divergence 
\[
    \widehat{C}_B(u,v) = \sum_{i=1}^r\sum_{k=1}^{\infty}
    \frac{(u_i^k-v_i^k)^2}{2k}=\sum_{i=1}^r
    \log\cosh\bigl(\operatorname{artanh}u_i-\operatorname{artanh}v_i\bigr) \qquad \text{for all } u, v \in \rho(\Omega).
\]

\subsubsection{Calabi's diastasis in Fermionic pairing state}
A more intricate case arise when we consider fermionic Bogoliubov and
Hartree--Fock--Bogoliubov states \cite{Robledo2009,CarlssonRotureau2021}. Let $\p=\Alt(\C^{2n})$ and $\Omega=\Alt(\C^{2n})_{<1}$. A matrix $X \in \Alt(\C^{2n})$ parametrizes normalized pure pairing state by
\[
\ket{\psi_X}
=\det(I+X^*X)^{-1/4}
\exp \Biggl(
\frac12\sum_{j,k=1}^{2n}
X_{jk}c_j^\dagger c_k^\dagger
\Biggr) \ket{0}.
\]
For $X, Y \in \Alt(\C^{2n})_{<1}$, the fidelity between the states is again measured by the inner product
\[
\lvert \braket{\psi_X|\psi_Y} \rvert^2
=
\frac{|\det(I+X^*Y)|}
{\sqrt{\det(I+X^*X)\det(I+Y^*Y)}}.
\]
Consider 
\[
    C_{\mathsf F}(X,Y)
    \coloneqq
    -\log\lvert\braket{\psi_X|\psi_Y} \rvert^2 \qquad \text{for all } X, Y \in \Alt(\C^{2n})_{<1},
\]
which is the Calabi diastasis induced by the potential
\[
\Phi_{\mathsf F}(Z)
=\frac12\log\det(I+Z^*Z).
\]
By Cauchy--Schwarz inequality, $C_{\mathsf F}(X,Y)\geq0$, with equality
if and only if $X=Y$. Moreover, $C_{\mathsf F}$ is smooth on $\Omega \times \Omega$ as $\lVert X^*Y\rVert_{\mathsf{op}}<1$. For every nonzero $H\in\p$,
\[
\frac{\partial^2}{\partial s\,\partial\overline{s}}
\Phi_{\mathsf F}(X+sH)
\Big\vert_{s=0}
=
\frac12\Bigl \lVert
(I+XX^*)^{-1/2}H(I+X^*X)^{-1/2}
\Bigr \rVert_{\mathsf F}^2>0.
\]
Therefore, $\Phi_{\mathsf F}$ defines a positive-definite K\"ahler metric and
$C_{\mathsf F}$ is a divergence. By unitary invariance of trace, $C_{\mathsf{F}}$ is $\K$-invariant. Furthermore, applying \cite[Theorem~3.3.13(c) and Theorem~3.3.14(c)]{HJ1991} to $f(x) = \log(1+x)$ yields
\[
\lvert \det(I+X^*Y) \rvert
\leq
\prod_{j=1}^{2n}
\bigl(1+\sigma_j(X)\sigma_j(Y)\bigr).
\]
Therefore, $C_{\mathsf F}(X,Y) \geq C_{\mathsf F}(X_+,Y_+)$ and $C_{\mathsf F}$ is $\K$-isotropic by Proposition~\ref{prop:char} with reduced divergence 
\[
\widehat C_{\mathsf F}(u,v)
=-2\sum_{i=1}^n
\log\cos\bigl(
\arctan u_i-\arctan v_i
\bigr) \qquad \text{for all } u, v \in \rho(\Omega).
\]

\subsection{Other isotropic divergences}

There are yet other isotropic divergences arising in different applications. We introduce three of them over the domain $\Omega=\Herm_+(\F^n)$.

\subsubsection{Bures--Wasserstein distance}

The squared Bures--Wasserstein distance \cite{BhatiaJainLim2019} is defined by
\[
D_{\mathsf{BW}}(X,Y)
=\tr X+\tr Y
-2\tr\bigl[(X^{1/2}YX^{1/2})^{1/2}\bigr]
\qquad\text{for all }X,Y\in\Herm_+(\F^n).
\]
On the positive-definite cone, this squared distance is a smooth divergence with a positive-definite induced Riemannian metric \cite[Section~4]{BhatiaJainLim2019}. By \cite[Theorem~3.3.14(a)]{HJ1991},
\[
\tr\bigl[(X^{1/2}YX^{1/2})^{1/2}\bigr]
=\lVert X^{1/2}Y^{1/2}\rVert_*
\leq\sum_{i=1}^n
\sqrt{\lambda_i(X)\lambda_i(Y)},
\]
so $D_{\mathsf{BW}}(X,Y)\geq
D_{\mathsf{BW}}(X_+,Y_+)$.
Thus $D_{\mathsf{BW}}$ is $\K$-isotropic with reduced divergence
\[
\widehat D_{\mathsf{BW}}(u,v)
=\sum_{i=1}^n(\sqrt{u_i}-\sqrt{v_i})^2
\qquad\text{for all }u,v\in\rho(\Omega).
\]

\subsubsection{Quantum Hellinger divergence}

The quantum Hellinger divergence \cite{LuoZhang2004,BhatiaGaubertJain2019} is defined by
\[
D_{\mathsf{QH}}(X,Y)
=\lVert X^{1/2}-Y^{1/2}\rVert_{\mathsf F}^2.
\]
Since the matrix square root is a smooth diffeomorphism of $\Omega$, $D_{\mathsf{QH}}$ is a divergence with a positive-definite induced metric. It is $\K$-invariant by unitary invariance of the Frobenius norm. By \cite[Theorem~3.3.14(a)]{HJ1991},
\[
\tr(X^{1/2}Y^{1/2})
\leq\sum_{i=1}^n
\sqrt{\lambda_i(X)\lambda_i(Y)},
\]
which implies $D_{\mathsf{QH}}(X,Y)\geq D_{\mathsf{QH}}(X_+,Y_+)$. Thus $D_{\mathsf{QH}}$ is $\K$-isotropic. Intriguingly, a direct computation shows that its reduced divergence coincides with the Bures--Wasserstein distance:
\[
\widehat D_{\mathsf{QH}}(u,v)
=\widehat D_{\mathsf{BW}}(u,v)
=\sum_{i=1}^n(\sqrt{u_i}-\sqrt{v_i})^2
\qquad\text{for all }u,v\in\rho(\Omega).
\]

\subsubsection{Quantum $q$-divergence}

Lastly, for $0<q<1$, consider the quantum $q$-divergence \cite{FeliceAy2019,HiaiEtAl2011} defined by
\[
D_{\mathsf{Q},q}(X,Y)
=\frac{1}{1-q}\tr\bigl(
qX+(1-q)Y-X^qY^{1-q}\bigr)
\quad\text{for }X,Y\in\Herm_+(\F^n).
\]
By \cite[Theorems~3.3.13(a') and 3.3.14(a)]{HJ1991} and H\"older's inequality 
\[
\tr(X^qY^{1-q})
\leq(\tr X)^q(\tr Y)^{1-q}
\leq q\tr X+(1-q)\tr Y,
\]
with equality throughout if and only if $X=Y$. Moreover, a direct calculation shows $D_{\mathsf{Q},q}$ is smooth with a positive-definite induced metric on $\Omega$. Thus $D_{\mathsf{Q},q}$ is a divergence.

If we restrict ourselves to the subset of density matrices, then the quantum $q$-divergence is the quantum Tsallis relative entropy \cite{Abe2003,FuruichiYanagiKuriyama2004}. Notice that it relates to the quantum Hellinger
and von Neumann divergences by
\[
D_{\mathsf{Q},1/2}=D_{\mathsf{QH}},
\qquad
\lim_{q\uparrow1}D_{\mathsf{Q},q}
=D_{\mathsf{vN}}.
\]
By \cite[Theorem~3.3.14(a)]{HJ1991}, for $X,Y\in\Herm_+(\F^n)$,
\[
\tr(X^qY^{1-q})
\leq\sum_{i=1}^n
\lambda_i(X)^q\lambda_i(Y)^{1-q},
\]
so $D_{\mathsf{Q},q}(X,Y)\geq D_{\mathsf{Q},q}(X_+,Y_+)$. Hence, the quantum $q$-divergence is $\K$-isotropic, with reduced divergence
\[
\widehat D_{\mathsf{Q},q}(u,v)
=\frac{1}{1-q}\sum_{i=1}^n
\bigl(qu_i+(1-q)v_i-u_i^qv_i^{1-q}\bigr)
\qquad\text{for all }u,v\in\rho(\Omega).
\]

\subsection{Separability of reduced divergences}

Notice that all reduced divergences in Table~\ref{tab:isotropic-divergences} are
\emph{separable}, i.e., for each divergence $D$, there exist an interval $\mathcal I_D\subseteq\R$ containing $\rho_i(A)$ for every $A\in\Omega$ and $i=1,\dots,r$, and a function $d_D:\mathcal I_D\times\mathcal I_D\to[0,\infty)$ satisfying $d_D(x,x)=0$ for all $x\in\mathcal I_D$, such that for all $u,v\in\rho(\Omega)$, $\widehat D(u,v)=d_D(u_1,v_1)+ \dots + d_D(u_r, v_r)$. As we will repeatedly use such a decomposition in the calculations in Section~\ref{sec:solve}, we summarize the corresponding $d_D$'s in Table~\ref{tab:d}.

\begin{table}[bth]
\centering
\footnotesize
\tabulinesep=0.6ex
\renewcommand{\arraystretch}{1.15}
\begin{adjustbox}{max width=\textwidth}
\begin{tabu}{@{}c|c|c@{}}
$D$ & $d_D(s,t)$ & $d_D(t,s)$\\
\hline

$D_{\mathsf F}$
& \multicolumn{2}{c}{$\tfrac{c_\p}{2}(s-t)^2$}\\

$D_{\mathsf{vN}}$
& $s\log\frac{s}{t}-s+t$
& $t\log\frac{t}{s}-t+s$\\

$D_{\mathsf{LD}}$
& $\frac{s}{t}-\log\frac{s}{t}-1$
& $\frac{t}{s}-\log\frac{t}{s}-1$\\

$D_{\mathsf{GKL}}$
& $\tfrac12\log\frac{1-t^2}{1-s^2}-\frac{t(s-t)}{1-t^2}$
& $\tfrac12\log\frac{1-s^2}{1-t^2}-\frac{s(t-s)}{1-s^2}$\\

$D_\beta$
& $f_\beta(s)-f_\beta(t)-(s-t)\operatorname{arsinh}(t/\beta)$
& $f_\beta(t)-f_\beta(s)-(t-s)\operatorname{arsinh}(s/\beta)$\\

\hline

$D_{\mathsf S}$
& \multicolumn{2}{c}{$\log\frac{s+t}{2\sqrt{st}}$}\\

$D_{\mathsf{QJS}}$
& \multicolumn{2}{c}{$\tfrac12\bigl[
s\log\frac{2s}{s+t}
+t\log\frac{2t}{s+t}\bigr]$}\\

\hline

$C_f$
& \multicolumn{2}{c}{$c_\p\displaystyle\sum_{k=1}^\infty
a_k(s^k-t^k)^2$}\\

$C_{\mathsf B}$
& \multicolumn{2}{c}{$\log\cosh(
\operatorname{artanh}s-\operatorname{artanh}t)$}\\

$C_{\mathsf F}$
& \multicolumn{2}{c}{$-2\log\cos(
\arctan s-\arctan t)$}\\

\hline

$D_{\mathsf{BW}}$
& \multicolumn{2}{c}{$(\sqrt{s}-\sqrt{t})^2$}\\

$D_{\mathsf{QH}}$
& \multicolumn{2}{c}{$(\sqrt{s}-\sqrt{t})^2$}\\

$D_{\mathsf{Q},q}$
& $\frac{qs+(1-q)t-s^qt^{1-q}}{1-q}$
& $\frac{qt+(1-q)s-t^qs^{1-q}}{1-q}$\\

\end{tabu}
\end{adjustbox}
\caption{Scalar functions of the separable reduced divergences.}
\label{tab:d}
\end{table}

\begin{lemma}\label{lem:separability}
For every $d_D$ in Table~\ref{tab:d}, let $f(x)=d_D(a,x)$ or $f(x)=d_D(x,a)$ for $x,a \in \mathcal{I}_D$. Then $f$ is strictly decreasing for $x<a$ and strictly
increasing for $x>a$.
\end{lemma}
\begin{proof}
It follows from the following derivatives.
\begin{enumerate}[\normalfont(i)]
\item For the Bregman divergences,
\[
\begin{aligned}
\frac{d}{dx}d_{\mathsf F}(a,x)
&=\frac{d}{dx}d_{\mathsf F}(x,a)
=c_\p(x-a),\\
\frac{d}{dx}d_{\mathsf{vN}}(a,x)
&=\frac{x-a}{x},&
\frac{d}{dx}d_{\mathsf{vN}}(x,a)
&=\log\frac{x}{a},\\
\frac{d}{dx}d_{\mathsf{LD}}(a,x)
&=\frac{x-a}{x^2},&
\frac{d}{dx}d_{\mathsf{LD}}(x,a)
&=\frac{x-a}{ax},\\
\frac{d}{dx}d_{\mathsf{GKL}}(a,x)
&=\frac{(x-a)(1+x^2)}{(1-x^2)^2},&
\frac{d}{dx}d_{\mathsf{GKL}}(x,a)
&=\frac{x}{1-x^2}-\frac{a}{1-a^2},\\
\frac{d}{dx}d_\beta(a,x)
&=\frac{x-a}{\sqrt{x^2+\beta^2}},&
\frac{d}{dx}d_\beta(x,a)
&=\operatorname{arsinh}\frac{x}{\beta}
-\operatorname{arsinh}\frac{a}{\beta}.
\end{aligned}
\]

\item For the Jensen divergences,
\[
\frac{d}{dx}d_{\mathsf S}(a,x)
=\frac{d}{dx}d_{\mathsf S}(x,a)
=\frac{x-a}{2x(a+x)}, \qquad
\frac{d}{dx}d_{\mathsf{QJS}}(a,x)
=\frac{d}{dx}d_{\mathsf{QJS}}(x,a)
=\frac12\log\frac{2x}{a+x}.
\]

\item For the Calabi diastases,
\[
\begin{aligned}
\frac{d}{dx}d_{C_f}(a,x)
=\frac{d}{dx}d_{C_f}(x,a)
&=2c_\p\sum_{k=1}^\infty
ka_kx^{k-1}(x^k-a^k),\\
\frac{d}{dx}d_{C_{\mathsf B}}(a,x)
=\frac{d}{dx}d_{C_{\mathsf B}}(x,a)
&=\frac{x-a}{(1-ax)(1-x^2)},\\
\frac{d}{dx}d_{C_{\mathsf F}}(a,x)
=\frac{d}{dx}d_{C_{\mathsf F}}(x,a)
&=\frac{2(x-a)}{(1+ax)(1+x^2)}.
\end{aligned}
\]

\item For the remaining divergences,
\[
\begin{aligned}
\frac{d}{dx}d_{\mathsf{BW}}(a,x)
&=\frac{d}{dx}d_{\mathsf{BW}}(x,a)
=1-\sqrt{\frac ax},
&\frac{d}{dx}d_{\mathsf{QH}}(a,x)
&=\frac{d}{dx}d_{\mathsf{QH}}(x,a)
=1-\sqrt{\frac ax},\\
\frac{d}{dx}d_{\mathsf{Q},q}(a,x)
&=1-\left(\frac ax\right)^q,&
\frac{d}{dx}d_{\mathsf{Q},q}(x,a)
&=\frac{q}{1-q}
\left[1-\left(\frac ax\right)^{1-q}\right]. 
\end{aligned}\qedhere
\]
\end{enumerate}
\end{proof}

\section{Isotropic nearness problems}\label{sec:solve}

We are now ready to solve a range of matrix nearness problems with isotropic divergences in Section~\ref{sec:divergences} and various invariant constraints. We solve them using the following three steps:
\begin{enumerate}[\normalfont(i)]
    \item Reduce the matrix nearness problems to vector nearness problems
    \begin{equation}\label{eq:reduced}
    \min_{v \in \rho(\mathcal{S})} \widehat{D}(\rho(A),v) \quad \text{or} \quad \min_{v \in \rho(\mathcal{S})} \widehat{D}(v,\rho(A))
    \end{equation}
    using Theorem~\ref{thm:reduce};
    \item solve the vector problem \eqref{eq:reduced} in closed form or by solve scalar equations;
    \item Recover the solution to \eqref{eq:near} using Corollary~\ref{cor:recover}.
\end{enumerate}
For the rest of this section, we denote $\lambda_i \coloneqq \lambda_i(A)$ if $A \in \Herm(\F^n)$ and $\sigma_i \coloneqq \sigma_i(A)$ in the other cases for $i = 1,\dots,n$ to avoid cluttering.

\subsection{rank constrained approximation}
The first nearness problem we consider is low-rank approximation over $\p = \F^{m \times n}$ and $\Omega \subseteq \F^{m \times n}$ open and $\K$-invariant. For $r = 1,\dots,n$, let 
\[
\mathcal{S}_r = \{X \in \Omega : \rank(X) \leq r\}.
\]
For $A \in \F^{m \times n}$, let $A = U \Sigma V^*$ be its singular value decomposition. For $r = 1,\dots,n$, let
\[
A_r \coloneqq U\diag_{m,n}(\sigma_1,\dots,\sigma_r,0,\dots,0)V^*.
\]
The following result can be viewed as a generalization of the Eckart--Young theorem \cite{Eck} to other isotropic divergences. The case of $D_{\mathsf{GKL}}(A,X)$ also follows from \cite[Theorem~2]{ChenPedersenAndersen2025}.

\begin{proposition}[Low-rank approximation]\label{prop:low-rank}
Let $\Omega\subseteq\F^{m\times n}$ be open and $\K$-invariant such that $A_k\in\Omega$ for every $A\in\Omega$ and $1\leq k\leq n$. Let $D$ be a $\K$-isotropic divergence on $\Omega$ with separable reduced divergence. Then for every $A \in \Omega$ and $1\leq r\leq n$,
\[
A_r\in\argmin_{X \in \mathcal{S}_r}D(A,X)
\quad\text{and}\quad
A_r\in\argmin_{X \in \mathcal{S}_r}D(X,A)
\]
with
\[
\min_{X \in \mathcal{S}_r}D(A,X)
=\sum_{i=r+1}^n d_D(\sigma_i,0),
\qquad
\min_{X \in \mathcal{S}_r}D(X,A)
=\sum_{i=r+1}^n d_D(0,\sigma_i).
\]
In particular, the minima for the squared Frobenius norm $D_\mathsf{F}$, $D_{\mathsf{GKL}}$, and Calabi's diastasis $C_f$ over the respective $\Omega$ are given by Table~\ref{tab:low-rank}.
\begin{table}[h]
\centering
\footnotesize
\tabulinesep=0.5ex
\begin{tabu}{@{}c|c|c@{}}
$D$
& $\min D(A,X)$
& $\min D(X,A)$\\
\hline

$D_{\mathsf F}$
& $\displaystyle\frac12\sum\limits_{i=r+1}^{n}\sigma_i^2$
& $\displaystyle\frac12\sum\limits_{i=r+1}^{n}\sigma_i^2$\\

$D_{\mathsf{GKL}}$
& $\displaystyle\frac12\sum\limits_{i=r+1}^{n}
\log\frac1{1-\sigma_i^2}$
& $\displaystyle\frac12\sum\limits_{i=r+1}^{n}
\left[\log(1-\sigma_i^2)
+\frac{2\sigma_i^2}{1-\sigma_i^2}\right]$\\

$C_f$
& $\displaystyle\sum\limits_{i=r+1}^{n}f(\sigma_i^2)$
& $\displaystyle\sum\limits_{i=r+1}^{n}f(\sigma_i^2)$

\end{tabu}
\caption{Rank-$r$ approximation errors.}
\label{tab:low-rank}
\end{table}
    
\end{proposition}
\begin{proof}
Since $\rank(X) \leq r$ if and only if $\sigma_{i}(X) = 0$ for all $i >r$,
\[
\sigma(\mathcal{S}_r)
=\{v\in\sigma(\Omega):v_{r+1}=\cdots=v_n=0\}.
\]
Therefore, for every $v\in\sigma(\mathcal{S}_r)$, 
\[
\widehat D(\sigma(A),v)
=\sum_{i=1}^r d_D(\sigma_i,v_i)
+\sum_{i=r+1}^n d_D(\sigma_i,0)
\geq\sum_{i=r+1}^n d_D(\sigma_i,0),
\]
with equality attained at $v_* \coloneqq (\sigma_1,\dots,\sigma_r,0,\dots,0)^\tp$. By assumption, $v_*=\sigma(A_r)\in\sigma(\mathcal{S}_r)$. Similarly, $\widehat D(v,\sigma(A)) \geq d_D(0,\sigma_{r+1})+ \dots + d_D(0,\sigma_n)$ with equality at $v_*$. Applying Theorem~\ref{thm:reduce} and Corollary~\ref{cor:recover} completes the proof.
\end{proof}

\subsection{Spectral constraints}

For $0\leq a\leq b$, let
\[
\mathcal S_{[a,b]}\coloneqq
\{X\in\Omega:\mathsf{\Lambda}(X)\cap[a,b]\neq\varnothing\},
\]
where $\p=\Herm(\F^n)$ and $\Omega=\Herm(\F^n)$ or $\Herm_+(\F^n)$. In the latter case, we assume $b>0$ to make sure $\mathcal{S}_{[a,b]} \neq \varnothing$. In particular, if $a=b$, the nearness problem over $\mathcal S_{\{a\}}$ is the prescribed eigenvalue problem, famously solved by Wilkinson under the Frobenius norm \cite{Wilkinson1986}.

\begin{proposition}[Spectral interval approximation]
\label{prop:spectral-interval}
Let $D$ be a $\K$-isotropic divergence defined over $\Omega=\Herm(\F^n)$ or $\Omega=\Herm_+(\F^n)$ with separable reduced divergence such that for every $t\in\mathcal I_D$, both $d_D(t,\cdot)$ and $d_D(\cdot,t)$ are strictly decreasing below $t$ and strictly increasing above $t$ on $\mathcal I_D$. If $A\in\Omega$ has eigenvalue decomposition $A=U\Lambda U^*$ such that $\mathsf{\Lambda}(A)\cap[a,b]=\varnothing$, then for $p\coloneqq\lvert \{i:\lambda_i>b\} \rvert$ and
\[
\begin{aligned}
A_b&\coloneqq
U\diag(\lambda_1,\dots,\lambda_{p-1},b,
\lambda_{p+1},\dots,\lambda_n)U^*,\\
A_a&\coloneqq
U\diag(\lambda_1,\dots,\lambda_p,a,
\lambda_{p+2},\dots,\lambda_n)U^*,
\end{aligned}
\]
whenever they are well-defined, we have 
\begin{align}
\min_{X\in\mathcal S_{[a,b]}}D(A,X)
&=\min\bigl\{d_D(\lambda_p,b),\,
  d_D(\lambda_{p+1},a)\bigr\}, \label{eq:spectral-interval-min} \\
\min_{X\in\mathcal S_{[a,b]}}D(X,A)
&=\min\bigl\{d_D(b,\lambda_p),\,
  d_D(a,\lambda_{p+1})\bigr\}, \label{eq:spectral-interval-min2}
\end{align}
where the $b$- and $a$-terms are included only when $p>0$ and $p<n$, respectively. When $0<p<n$, $A_b$ is a minimizer if $d_D(\lambda_p,b)\leq d_D(\lambda_{p+1},a)$, and $A_a$ otherwise. The case of \eqref{eq:spectral-interval-min2} follows analogously, with the arguments of $d_D$ reversed.
\end{proposition}
\begin{proof}

By Theorem~\ref{thm:reduce}, it suffices to solve the corresponding reduced nearness problems. Let $v\in\lambda(\mathcal S_{[a,b]})$. Then $v_j \in [a,b]$ for some $j$. If $\lambda_j>b$, then
\[
\widehat D(\lambda(A),v)
\geq d_D(\lambda_j,v_j)
\geq d_D(\lambda_j,b)
\geq d_D(\lambda_p,b)
\]
where the last two inequalities follow from the assumption. If $\lambda_j<a$, the same argument gives that $\widehat D(\lambda(A),v) \geq d_D(\lambda_{p+1},a)$. Since the lower bound is attained by the eigenvalues of $A_b$ or $A_a$, we obtain \eqref{eq:spectral-interval-min}. The equality \eqref{eq:spectral-interval-min2} follows from the same argument. The minimizers are obtained by Corollary~\ref{cor:recover}. 
\end{proof}

\subsection{Condition-number constraint}

Let $\p=\F^{m\times n}$, $\kappa\geq1$, and $\Omega=\F^{m \times n}$ for $D_{\mathsf F}$, $\Omega=\R^{m\times n}_{<1}$ for $D_{\mathsf{GKL}}$, and $\Omega=\C^{m\times n}_{<\mu}$ for $C_f$. Consider the condition-number-constrained approximation
\begin{equation}\label{eq:cond}
    \min_{X \in \mathcal{S}_\kappa} D(A, X), \quad \min_{X \in \mathcal{S}_\kappa} D(X, A), \quad \mathcal{S}_\kappa \coloneqq \{X \in \Omega : \sigma_1(X) \leq \kappa \sigma_n(X)\}.
\end{equation} 
\begin{proposition}\label{prop:cond}
Let $D= D_\mathsf{F}$, $D_{\mathsf{GKL}}$, or $C_f$, with $\Omega$ as given above. Let $A \in \Omega \setminus \mathcal{S}_\kappa$ with singular value decomposition $A = U \Sigma V^*$. For $t \geq 0$, let $\mathcal{I}(t) \coloneqq \{i: \sigma_i > \kappa t\}$, $\mathcal{J}(t) \coloneqq \{i:\sigma_i < t\}$, and define $v(t)$ by $v_i(t)\coloneqq \min\bigl\{\kappa t, \max\{\sigma_i,t\} \bigr\}$. 

If $t_* \in (\sigma_n, \sigma_1/\kappa)$ is a solution to the corresponding equation in Table~\ref{tab:cond-right} that attains the minimal $\widehat{D}\bigl(\sigma(A),v(t)\bigr)$ among all solutions, then $X_* \coloneqq U\diag_{m,n}\bigl(v(t_*)\bigr) V^*$ is a minimizer to \eqref{eq:cond} on $D(A,X)$ with minimum given in Table~\ref{tab:cond-right} at $v=v(t_*)$. A minimizer to \eqref{eq:cond} on $D(X,A)$ is similarly obtained with Table~\ref{tab:cond-left}.

\begin{table}[htb]
\centering\scriptsize
\tabulinesep=0.7ex
\renewcommand{\arraystretch}{1.13}
\begin{adjustbox}{max width=\textwidth}
\begin{tabu}{@{}c|c|c@{}}
$D$ & $t$ & $\min D(A,X)$\\\hline

$D_{\mathsf F}$
& $\displaystyle
t=\frac{\kappa\sum_{i\in\mathcal I(t)}\sigma_i+
\sum_{i\in\mathcal J(t)}\sigma_i}
{\kappa^2\lvert\mathcal I(t)\rvert+
\lvert\mathcal J(t)\rvert}$
& $\displaystyle
\frac12\sum_{i=1}^n(\sigma_i-v_i)^2$\\

$D_{\mathsf{GKL}}$
& $\displaystyle
\kappa\left(\kappa t\lvert\mathcal I(t)\rvert-
\sum_{i\in\mathcal I(t)}\sigma_i\right)
\frac{1+\kappa^2t^2}{(1-\kappa^2t^2)^2}
+\left(t\lvert\mathcal J(t)\rvert-
\sum_{i\in\mathcal J(t)}\sigma_i\right)
\frac{1+t^2}{(1-t^2)^2}=0$
& $\displaystyle
\frac12\sum_{i=1}^n\biggl[
\log\frac{1-v_i^2}{1-\sigma_i^2}
-\frac{2v_i(\sigma_i-v_i)}
{1-v_i^2}
\biggr]$\\

$C_f$
& $\displaystyle
\sum_{j=1}^\infty ja_jt^{j-1}
\biggl[
\left(\kappa^{2j}\lvert\mathcal I(t)\rvert+
\lvert\mathcal J(t)\rvert\right)t^j
-\kappa^j\sum_{i\in\mathcal I(t)}\sigma_i^j
-\sum_{i\in\mathcal J(t)}\sigma_i^j
\biggr]=0$
& $\displaystyle
\sum_{i=1}^n\sum_{j=1}^\infty a_j
(\sigma_i^j-v_i^j)^2$\\

\end{tabu}
\end{adjustbox}
\caption{Solution to the right condition number-constrained nearness problem.}
\label{tab:cond-right}
\end{table}

\begin{table}[htb]
\centering\scriptsize
\tabulinesep=0.7ex
\renewcommand{\arraystretch}{1.13}
\begin{adjustbox}{max width=\textwidth}
\begin{tabu}{@{}c|c|c@{}}
$D$ & $t$ & $\min D(X,A)$\\\hline

$D_{\mathsf F}$
& $\displaystyle
t=\frac{\kappa\sum_{i\in\mathcal I(t)}\sigma_i+
\sum_{i\in\mathcal J(t)}\sigma_i}
{\kappa^2\lvert\mathcal I(t)\rvert+
\lvert\mathcal J(t)\rvert}$
& $\displaystyle
\frac12\sum_{i=1}^n(\sigma_i-v_i)^2$\\

$D_{\mathsf{GKL}}$
& $\displaystyle
\kappa\left(
\frac{\kappa t\lvert\mathcal I(t)\rvert}
{1-\kappa^2t^2}
-\sum_{i\in\mathcal I(t)}
\frac{\sigma_i}{1-\sigma_i^2}\right)
+\frac{t\lvert\mathcal J(t)\rvert}{1-t^2}
-\sum_{i\in\mathcal J(t)}
\frac{\sigma_i}{1-\sigma_i^2}=0$
& $\displaystyle
\frac12\sum_{i=1}^n\biggl[
\log\frac{1-\sigma_i^2}{1-v_i^2}
-\frac{2\sigma_i(v_i-\sigma_i)}
{1-\sigma_i^2}
\biggr]$\\

$C_f$
& $\displaystyle
\sum_{j=1}^\infty ja_jt^{j-1}
\biggl[
\left(\kappa^{2j}\lvert\mathcal I(t)\rvert+
\lvert\mathcal J(t)\rvert\right)t^j
-\kappa^j\sum_{i\in\mathcal I(t)}\sigma_i^j
-\sum_{i\in\mathcal J(t)}\sigma_i^j
\biggr]=0$
& $\displaystyle
\sum_{i=1}^n\sum_{j=1}^\infty a_j
(\sigma_i^j-v_i^j)^2$\\

\end{tabu}
\end{adjustbox}
\caption{Solution to the left condition number-constrained nearness problem.}
\label{tab:cond-left}
\end{table}
\end{proposition}  
\begin{proof}
By Theorem~\ref{thm:reduce}, it suffices to minimize the corresponding reduced divergences over 
\[
\sigma(\mathcal S_\kappa)
=\{v\in\sigma(\Omega):v_1\leq\kappa v_n\}
=\bigcup_{t\geq0}
\bigl(\sigma(\Omega)\cap[t,\kappa t]^n\bigr),
\]
where the last equality follows from the following observation: If $v\in\sigma(\mathcal{S}_\kappa)$, then $v\in[v_n,\kappa v_n]^n$. Conversely, if $v\in\sigma(\Omega)\cap[t,\kappa t]^n$ for some $t\geq0$, then $v_1\leq\kappa t\leq\kappa v_n$. Therefore, we first minimize the reduced
divergences over $\sigma(\Omega)\cap[t,\kappa t]^n$ for feasible $t\geq0$, and then choose the optimal $t$.

For fixed $t\geq0$ such that $\sigma(\Omega)\cap[t,\kappa t]^n\neq\varnothing$, the functions $d_D(\sigma_i,x)$ and $d_D(x,\sigma_i)$ are minimized over $x\in[t,\kappa t]$ at $x=v_i(t)$ by Lemma~\ref{lem:separability}. Since $\sigma_1\geq\cdots\geq\sigma_n$,
we have
\[
t\leq v_n(t)\leq\cdots\leq v_1(t)\leq\kappa t.
\]
Moreover, whenever $\sigma(\Omega)\cap[t,\kappa t]^n$ is nonempty, $v(t)\in\sigma(\Omega)$ for the domains considered here. Therefore, 
\[
\min_{v\in\sigma(\Omega)\cap[t,\kappa t]^n}
\widehat D(\sigma(A),v)
=\sum_{i=1}^n d_D(\sigma_i,v_i(t)), \quad
\min_{v\in\sigma(\Omega)\cap[t,\kappa t]^n}
\widehat D(v,\sigma(A))
=\sum_{i=1}^n d_D(v_i(t),\sigma_i).
\]

If $t<\sigma_n$, then $v_i(t)\leq v_i(\sigma_n)\leq\sigma_i$ for every $i = 1,\dots,n$. Thus, by
Lemma~\ref{lem:separability}, replacing $t$ by $\sigma_n$ does not increase either reduced divergence. Similarly, if $t>\sigma_1/\kappa$, then $\sigma_i\leq v_i(\sigma_1/\kappa)\leq v_i(t)$, so replacing $t$ by $\sigma_1/\kappa$ does not increase either reduced divergence. Therefore, 
\[
\min_{X\in\mathcal S_\kappa}D(A,X)
=\min_{\sigma_n\leq t\leq\sigma_1/\kappa}
\sum_{i=1}^n d(\sigma_i,v_i(t)), \quad
\min_{X\in\mathcal S_\kappa}D(X,A)
=\min_{\sigma_n\leq t\leq\sigma_1/\kappa}
\sum_{i=1}^n d(v_i(t),\sigma_i).
\]

Recall that $\sigma_1>\kappa\sigma_n$ since $A\notin\mathcal S_\kappa$. Notice that the two reduced objective functions are continuous on the compact interval$[\sigma_n,\sigma_1/\kappa]$, so both attain their minima. Moreover, by the derivatives in the proof of Lemma~\ref{lem:separability}, both objective functions have strictly negative right derivatives at $t=\sigma_n$ and strictly positive left derivatives at $t=\sigma_1/\kappa$. Thus, both minima occur in the open interval $(\sigma_n,\sigma_1/\kappa)$. Differentiating $\sum_{i=1}^n d(v_i(t),\sigma_i)$ and $\sum_{i=1}^n d(\sigma_i, v_i(t))$ with respect to $t$ yields the equations in Table~\ref{tab:cond-right} and Table~\ref{tab:cond-left}. If an equation admits more than one solution, choosing one that minimizes the corresponding reduced divergence among all solutions yields a global minimizer. Lastly, we obtain $X_*$ by Corollary~\ref{cor:recover}.
\end{proof}
We would like to remark that the Frobenius case of Proposition~\ref{prop:cond} can also be obtained from \cite{Nievergelt}.

\subsection{Trace constraint}

Let $\p=\Herm(\F^n)$, $A\in\Herm_+(\F^n)$, and $\tau>0$. We consider $D=D_{\mathsf F}$, $D_{\mathsf{vN}}$, $D_{\mathsf{LD}}$, $D_\beta$, $D_{\mathsf S}$, $D_{\mathsf{QJS}}$, $D_{\mathsf{BW}}$, $D_{\mathsf{QH}}$, or $D_{\mathsf{Q},q}$, each restricted to $\Omega=\Herm_+(\F^n)$. We may extend these
divergences to positive semidefinite matrices by 
\[
D(A,X)\coloneqq
\lim_{\varepsilon\downarrow0}D(A,X+\varepsilon I),
\qquad
D(X,A)\coloneqq
\lim_{\varepsilon\downarrow0}D(X+\varepsilon I,A),
\]
allowing the value to be $+\infty$. Consider the nearness problems
\begin{equation}\label{eq:trace}
    \min_{X\in\mathcal S_\tau}D(A,X),\quad
    \min_{X\in\mathcal S_\tau}D(X,A),\quad
    \mathcal S_\tau\coloneqq
    \{X\in\p:X\succeq0,\,\tr X=\tau\}.
\end{equation}
In particular, if $\F=\C$ and $\tau=1$, $\mathcal S_1$ is exactly the set of density matrices. The reduced problem in the following result has been studied in the cases of Euclidean distance and vector Bregman divergences \cite{Duchi2008,Krichene2015} and the case of approximating density matrices in Frobenius norm has been solved in \cite{Smolin2012}.

\begin{proposition}\label{prop:trace}
Let $D$ be one of the divergences listed above. Let $A\in\Herm_+(\F^n)$ have eigenvalue decomposition $A=U\Lambda U^*$ and $\tr(A)\neq\tau$. Let $x_i$ be given by Table~\ref{tab:trace-right} and Table~\ref{tab:trace-left}, where \begin{enumerate}[\normalfont(i)]
\item $\eta\in\R$ satisfies $\sum_{i=1}^n x_i=\tau$; \label{item:trace1}
\item for $D_{\mathsf{LD}}(A,X)$ and $D_{\mathsf S}$, the minus sign is taken for $i=2,\dots,n$, while either sign is allowed for $i=1$. Among all real, positive candidates satisfying \ref{item:trace1}, including all admissible multipliers $\eta$, choose one with the smallest corresponding objective value. \label{item:trace2}
\end{enumerate}
Then $X_*=U\diag(x_1,\dots,x_n)U^*$ is a minimizer of the corresponding problem in \eqref{eq:trace}, with minimum given in the same table.

\begin{table}[htb]
\centering\scriptsize
\tabulinesep=0.7ex
\renewcommand{\arraystretch}{1.13}
\begin{adjustbox}{max width=\textwidth}
\begin{tabu}{@{}c|c|c@{}}
$D$ & $x_i$ & $\min D(A,X)$\\
\hline

$D_{\mathsf F}$
& $\displaystyle\max\{\lambda_i+\eta,0\}$
& $\displaystyle\frac12\sum_{i=1}^n
(\lambda_i-x_i)^2$\\

$D_{\mathsf{vN}}$
& $\displaystyle\frac{\tau}{\tr(A)}\lambda_i$
& $\displaystyle\tr(A)\log\frac{\tr(A)}{\tau}
-\tr(A)+\tau$\\

$D_{\mathsf{LD}}$
& $\displaystyle
\frac{1\pm\sqrt{1-4\eta\lambda_i}}{2\eta}$
& $\displaystyle\sum_{i=1}^n\left[
\frac{\lambda_i}{x_i}
-\log\frac{\lambda_i}{x_i}-1
\right]$\\

$D_\beta$
& $\displaystyle\begin{cases}
0,&\eta\leq-\lambda_i/\beta,\\
\frac{\lambda_i^2-\beta^2}{2\lambda_i},
&\eta=-1,\ \lambda_i>\beta,\\
\frac{\lambda_i+\eta
\sqrt{\lambda_i^2+(1-\eta^2)\beta^2}}
{1-\eta^2},
&-\lambda_i/\beta<\eta<1,\quad\eta\neq-1
\end{cases}$
& $\displaystyle\sum_{i=1}^n\left[
f_\beta(\lambda_i)-f_\beta(x_i)
-(\lambda_i-x_i)
\operatorname{arsinh}\frac{x_i}{\beta}
\right]$\\

$D_{\mathsf S}$
& $\displaystyle
\frac{1-2\eta\lambda_i\pm
\sqrt{1-12\eta\lambda_i+4\eta^2\lambda_i^2}}
{4\eta}$
& $\displaystyle\sum_{i=1}^n
\log\frac{\lambda_i+x_i}
{2\sqrt{\lambda_i x_i}}$\\

$D_{\mathsf{QJS}}$
& $\displaystyle\frac{\tau}{\tr(A)}\lambda_i$
& $\displaystyle\frac12\left[
\tr(A)\log\frac{2\tr(A)}{\tr(A)+\tau}
+\tau\log\frac{2\tau}{\tr(A)+\tau}
\right]$\\

$D_{\mathsf{BW}}$
& $\displaystyle\frac{\tau}{\tr(A)}\lambda_i$
& $\displaystyle
\left(\sqrt{\tr(A)}-\sqrt\tau\right)^2$\\

$D_{\mathsf{QH}}$
& $\displaystyle\frac{\tau}{\tr(A)}\lambda_i$
& $\displaystyle
\left(\sqrt{\tr(A)}-\sqrt\tau\right)^2$\\

$D_{\mathsf{Q},q}$
& $\displaystyle\frac{\tau}{\tr(A)}\lambda_i$
& $\displaystyle\frac{
q\tr(A)+(1-q)\tau-\tr(A)^q\tau^{1-q}}
{1-q}$\\

\end{tabu}
\end{adjustbox}
\caption{Solution to the right trace-constrained nearness problem.}
\label{tab:trace-right}
\end{table}

\begin{table}[htb]
\centering\scriptsize
\tabulinesep=0.7ex
\renewcommand{\arraystretch}{1.13}
\begin{adjustbox}{max width=\textwidth}
\begin{tabu}{@{}c|c|c@{}}
$D$ & $x_i$ & $\min D(X,A)$\\
\hline

$D_{\mathsf F}$
& $\displaystyle\max\{\lambda_i+\eta,0\}$
& $\displaystyle\frac12\sum_{i=1}^n
(\lambda_i-x_i)^2$\\

$D_{\mathsf{vN}}$
& $\displaystyle\frac{\tau}{\tr(A)}\lambda_i$
& $\displaystyle\tau\log\frac{\tau}{\tr(A)}
-\tau+\tr(A)$\\

$D_{\mathsf{LD}}$
& $\displaystyle\frac{\lambda_i}
{1-\eta\lambda_i}$
& $\displaystyle\sum_{i=1}^n\left[
\frac{x_i}{\lambda_i}
-\log\frac{x_i}{\lambda_i}-1
\right]$\\

$D_\beta$
& $\displaystyle\max\left\{
\beta\sinh\left(
\operatorname{arsinh}\frac{\lambda_i}{\beta}
+\eta\right),0\right\}$
& $\displaystyle\sum_{i=1}^n\left[
f_\beta(x_i)-f_\beta(\lambda_i)
-(x_i-\lambda_i)
\operatorname{arsinh}\frac{\lambda_i}{\beta}
\right]$\\

$D_{\mathsf S}$
& $\displaystyle
\frac{1-2\eta\lambda_i\pm
\sqrt{1-12\eta\lambda_i+4\eta^2\lambda_i^2}}
{4\eta}$
& $\displaystyle\sum_{i=1}^n
\log\frac{\lambda_i+x_i}
{2\sqrt{\lambda_i x_i}}$\\

$D_{\mathsf{QJS}}$
& $\displaystyle\frac{\tau}{\tr(A)}\lambda_i$
& $\displaystyle\frac12\left[
\tr(A)\log\frac{2\tr(A)}{\tr(A)+\tau}
+\tau\log\frac{2\tau}{\tr(A)+\tau}
\right]$\\

$D_{\mathsf{BW}}$
& $\displaystyle\frac{\tau}{\tr(A)}\lambda_i$
& $\displaystyle
\left(\sqrt{\tr(A)}-\sqrt\tau\right)^2$\\

$D_{\mathsf{QH}}$
& $\displaystyle\frac{\tau}{\tr(A)}\lambda_i$
& $\displaystyle
\left(\sqrt{\tr(A)}-\sqrt\tau\right)^2$\\

$D_{\mathsf{Q},q}$
& $\displaystyle\frac{\tau}{\tr(A)}\lambda_i$
& $\displaystyle\frac{
q\tau+(1-q)\tr(A)-\tau^q\tr(A)^{1-q}}
{1-q}$\\

\end{tabu}
\end{adjustbox}
\caption{Solution to the left trace-constrained nearness problem.}
\label{tab:trace-left}
\end{table}
\end{proposition}

\begin{proof}
Notice that the extended objective functions are lower semicontinuous on the compact set $\mathcal S_\tau$, so both minima are attained. Let $X\in\mathcal S_\tau$ and $X_\varepsilon\coloneqq X + \varepsilon I\succ0$. By Proposition~\ref{prop:char}, for every $\varepsilon>0$,
\[
D(A,X_\varepsilon)\geq\widehat D\bigl(
\lambda(A),
\lambda(X)+\varepsilon(1,\dots,1)\bigr), \qquad
D(X_\varepsilon,A)
\geq\widehat D\bigl(
\lambda(X)+\varepsilon(1,\dots,1),
\lambda(A)\bigr).
\]
Passing to the limit $\varepsilon\downarrow0$, with the same boundary convention for $\widehat D$, gives
\[
D(A,X)\geq\widehat D(\lambda(A),\lambda(X)),
\qquad
D(X,A)\geq\widehat D(\lambda(X),\lambda(A)).
\]
On the other hand, equality is attained  at $U\diag(\lambda(X))U^*$. Therefore, we obtain the reduced problem
\begin{equation}\label{eq:trace-reduced}
\min_{X\in\mathcal S_\tau}D(A,X)
=\min_{\substack{
x_1\geq\cdots\geq x_n\geq0\\
\sum_{i=1}^n x_i=\tau}}
\widehat D(\lambda(A),x), \qquad
\min_{X\in\mathcal S_\tau}D(X,A)
=\min_{\substack{
x_1\geq\cdots\geq x_n\geq0\\
\sum_{i=1}^n x_i=\tau}}
\widehat D(x,\lambda(A)).
\end{equation}

Notice that for each corresponding scalar function $d_D$ in Table~\ref{tab:d}, $\partial_y\partial_xd_D(y,x)<0$ for $y,x>0$. Thus, if $a\geq b$ and $u\geq v$, then
\begin{equation}\label{eq:mix}
d_D(a,u)+d_D(b,v)
-d_D(a,v)-d_D(b,u)=
\int_b^a\int_v^u
\partial_y\partial_xd_D(y,x)
\,dx\,dy\leq0.
\end{equation}
Therefore, rearranging entries of $x$ in nonincreasing order does not increase the reduced divergence. The analogous statement holds with the arguments of $d_D$ reversed and extends
to zero coordinates. Therefore, we may first solve the scalar problems on the full simplex
$\{x\in\R^n:x_i\geq0,\, \sum_{i=1}^n x_i=\tau\}$ and then choose an ordering.

For $D_{\mathsf F}$, $D_{\mathsf{vN}}$, $D_\beta$, $D_{\mathsf{LD}}(X,A)$, $D_{\mathsf{QJS}}$, $D_{\mathsf{BW}}$, $D_{\mathsf{QH}}$, and $D_{\mathsf{Q},q}$, the corresponding scalar objectives are
strictly convex, as follows by differentiating the functions in Table~\ref{tab:d}. The KKT conditions at a minimizer are
\[
\begin{cases}
f_i'(x_i)=\eta,&x_i>0,\\
f_i'(0_+)\geq\eta,&x_i=0,
\end{cases}
\qquad
\sum_{i=1}^n x_i=\tau,
\]
where $f_i(x)=d_D(\lambda_i,x)$ or $f_i(x)=d_D(x,\lambda_i)$, depending on the orientation. In strictly convex cases, they are sufficient conditions for global optimality.

For $D_{\mathsf F}$, $D_\beta$, and $D_{\mathsf{LD}}(X,A)$, the KKT equations are
\[
\begin{aligned}
x_i-\lambda_i
&=\eta,
&&\text{for }D_{\mathsf F},\\
\frac{x_i-\lambda_i}
{\sqrt{x_i^2+\beta^2}}
&=\eta,
&&\text{for }D_\beta(A,X),\\
\operatorname{arsinh}\frac{x_i}{\beta}
-\operatorname{arsinh}\frac{\lambda_i}{\beta}
&=\eta,
&&\text{for }D_\beta(X,A),\\
\frac1{\lambda_i}-\frac1{x_i}
&=\eta,
&&\text{for }D_{\mathsf{LD}}(X,A).
\end{aligned}
\]
Solving these equations, together with the conditions at zero, gives the eigenvalues in Tables~\ref{tab:trace-right} and~\ref{tab:trace-left}.

For $D_{\mathsf{LD}}(A,X)$ and $D_{\mathsf S}$, the scalar objectives are not generally convex.
Nevertheless, they grow to $+\infty$ as any coordinate tends to zero, so every global minimizer lies in $\Herm_+(\F^n)$ and satisfies the KKT equations. Here $\eta\neq0$, since $\eta=0$ yields $x_i=\lambda_i$ for every $i$, contradicting $\tr(A)\neq\tau$. Solving these equations gives the quadratic roots in the tables. Their second derivatives are
\[
\frac{\partial^2}{\partial x^2}
d_{\mathsf{LD}}(\lambda,x)
=\frac{2\lambda-x}{x^3}, \qquad 
\frac{\partial^2}{\partial x^2}
d_{\mathsf S}(\lambda,x)
=\frac{\lambda^2+2\lambda x-x^2}
{2x^2(\lambda+x)^2}.
\]
Write $x_-(\lambda)$ and $x_+(\lambda)$ for the roots obtained with the minus and plus signs, respectively. Whenever these roots are distinct and positive, $f_\lambda''(x_+(\lambda))<0$. Thus, at most one coordinate of a constrained local minimizer can take a distinct plus root, since otherwise the second derivative along $e_i-e_j$ satisfies $f_i''(x_i)+f_j''(x_j)<0$. Moreover, $\partial_\lambda f_\lambda'(x)<0$ for both divergences. Implicit differentiation yields that
\[
\frac{d}{d\lambda}x_+(\lambda)
=-\frac{\partial_\lambda f_\lambda'(x_+(\lambda))}
{f_\lambda''(x_+(\lambda))}<0.
\]
Hence, if a distinct plus root occurred at $j>1$, then
\[
x_1=x_-(\lambda_1)
\leq x_+(\lambda_1)
\leq x_+(\lambda_j)=x_j,
\]
with at least one strict inequality, contradicting the nonincreasing ordering of the optimal coordinates. Therefore, only $i=1$ may take the plus root. Comparing all feasible candidates described in \ref{item:trace1} and \ref{item:trace2} yields a global minimizer.

For $D_{\mathsf{vN}}$, $D_{\mathsf{QJS}}$, $D_{\mathsf{BW}}$, $D_{\mathsf{QH}}$, and $D_{\mathsf{Q},q}$, the first derivatives in Lemma~\ref{lem:separability} are strictly increasing functions of
$x_i/\lambda_i$ and tend to $-\infty$ as $x_i\downarrow0$. Hence all optimal coordinates are positive, and the KKT equations imply that $x_i/\lambda_i$ is independent of $i$. Thus
\[
x_i=\frac{\tau}{\tr(A)}\lambda_i,
\qquad i=1,\dots,n.
\]
Applying Corollary~\ref{cor:recover} yields the minimizer.
\end{proof}

\subsection{Determinant constraint}

Let $\p=\Herm(\F^n)$, $\Omega=\Herm_+(\F^n)$, and $\delta>0$. Consider the nearness problems
\begin{equation}\label{eq:det}
    \min_{X\in\mathcal S_\delta}D(A,X),\qquad
    \min_{X\in\mathcal S_\delta}D(X,A),\qquad
    \mathcal S_\delta\coloneqq
    \{X\in\Herm_+(\F^n):\det X=\delta\},
\end{equation}
which arises in statistics and experimental design \cite{Frisch1929,Wilks1932,Wald1943}. The case of Frobenius norm has been considered in \cite{JS2026}.

\begin{proposition}\label{prop:determinant}
Let $D=D_{\mathsf F}$, $D_{\mathsf{vN}}$, $D_{\mathsf{LD}}$, $D_\beta$, $D_{\mathsf S}$, $D_{\mathsf{QJS}}$, $D_{\mathsf{BW}}$, $D_{\mathsf{QH}}$, or $D_{\mathsf{Q},q}$. Let $A\in\Herm_+(\F^n)$ satisfy
$\det A\neq\delta$, and let $A=U\Lambda U^*$ be its eigenvalue decomposition. If $x_1,\dots,x_n$ are defined according to Table~\ref{tab:det-right} or Table~\ref{tab:det-left} as detailed below, then $X_*=U\diag(x_1,\dots,x_n)U^*$ is a minimizer of the corresponding problem in \eqref{eq:det}, with minimum value given in the same table.

In Tables~\ref{tab:det-right} and~\ref{tab:det-left}, choose $\eta$ so that $\prod_{i=1}^n x_i=\delta$. If $\delta>\det A$, take the unique roots $x_i>\lambda_i$. If $\delta<\det A$, take
$0<x_i<\lambda_i$. When there are two distinct positive roots, take the larger root for
$i=1,\dots,n-1$ and compare both choices for $i=n$. For $D_{\mathsf{BW}}$ and $D_{\mathsf{QH}}$, only the plus sign is admissible when $\eta>0$. Among all admissible solutions of the product constraint, choose one with the smallest corresponding objective value. In the tables, $W_0$ and $W_{-1}$ denote the two real branches of the Lambert $W$ function.

\begin{table}[htb]
\centering\scriptsize
\tabulinesep=0.7ex
\renewcommand{\arraystretch}{1.13}
\begin{adjustbox}{max width=\textwidth}
\begin{tabu}{@{}c|c|c@{}}
$D$ &$x_i$ & $\min D(A,X)$\\
\hline

$D_{\mathsf F}$
& $\displaystyle
\frac{\lambda_i\pm\sqrt{\lambda_i^2+4\eta}}{2}$
& $\displaystyle
\frac12\sum_{i=1}^n(\lambda_i-x_i)^2$\\

$D_{\mathsf{vN}}$
& $\displaystyle\lambda_i+\eta$
& $\displaystyle\sum_{i=1}^n
\left[\lambda_i\log\frac{\lambda_i}{x_i}
-\lambda_i+x_i\right]$\\

$D_{\mathsf{LD}}$
& $\displaystyle
\left(\frac{\delta}{\det A}\right)^{1/n}\lambda_i$
& $\displaystyle n\left[
\left(\frac{\det A}{\delta}\right)^{1/n}
+\frac1n\log\frac{\delta}{\det A}-1
\right]$\\

$D_\beta$
& $\displaystyle
\frac{x_i(x_i-\lambda_i)}
{\sqrt{x_i^2+\beta^2}}=\eta$
& $\displaystyle\sum_{i=1}^n\left[
f_\beta(\lambda_i)-f_\beta(x_i)
-(\lambda_i-x_i)
\operatorname{arsinh}\frac{x_i}{\beta}
\right]$\\

$D_{\mathsf S}$
& $\displaystyle
\left(\frac{\delta}{\det A}\right)^{1/n}\lambda_i$
& $\displaystyle n\log
\frac{(\det A)^{1/n}+\delta^{1/n}}
{2(\delta\det A)^{1/(2n)}}$\\

$D_{\mathsf{QJS}}$
& $\displaystyle
\frac{x_i}{2}\log
\frac{2x_i}{\lambda_i+x_i}=\eta$
& $\displaystyle\frac12\sum_{i=1}^n\left[
\lambda_i\log\frac{2\lambda_i}{\lambda_i+x_i}
+x_i\log\frac{2x_i}{\lambda_i+x_i}
\right]$\\

$D_{\mathsf{BW}}$
& $\displaystyle
\frac{\bigl(\sqrt{\lambda_i}
\pm\sqrt{\lambda_i+4\eta}\bigr)^2}{4}$
& $\displaystyle
\sum_{i=1}^n(\sqrt{\lambda_i}-\sqrt{x_i})^2$\\

$D_{\mathsf{QH}}$
& $\displaystyle
\frac{\bigl(\sqrt{\lambda_i}
\pm\sqrt{\lambda_i+4\eta}\bigr)^2}{4}$
& $\displaystyle
\sum_{i=1}^n(\sqrt{\lambda_i}-\sqrt{x_i})^2$\\

$D_{\mathsf{Q},q}$
& $\displaystyle
x_i-\lambda_i^q x_i^{1-q}=\eta$
& $\displaystyle\frac1{1-q}\sum_{i=1}^n
\left[q\lambda_i+(1-q)x_i
-\lambda_i^q x_i^{1-q}\right]$\\

\end{tabu}
\end{adjustbox}
\caption{Solution to the right determinant-constrained nearness problem.}
\label{tab:det-right}
\end{table}

\begin{table}[htb]
\centering\scriptsize
\tabulinesep=0.7ex
\renewcommand{\arraystretch}{1.13}
\begin{adjustbox}{max width=\textwidth}
\begin{tabu}{@{}c|c|c@{}}
$D$ & $x_i$ & $\min D(X,A)$\\
\hline

$D_{\mathsf F}$
& $\displaystyle
\frac{\lambda_i\pm\sqrt{\lambda_i^2+4\eta}}{2}$
& $\displaystyle
\frac12\sum_{i=1}^n(\lambda_i-x_i)^2$\\

$D_{\mathsf{vN}}$
& $\displaystyle
\lambda_i\exp(W_0(\eta/\lambda_i))
\ \text{or}\
\lambda_i\exp(W_{-1}(\eta/\lambda_i))$
& $\displaystyle\sum_{i=1}^n\left[
x_i\log\frac{x_i}{\lambda_i}
-x_i+\lambda_i\right]$\\

$D_{\mathsf{LD}}$
& $\displaystyle
\left(\frac{\delta}{\det A}\right)^{1/n}\lambda_i$
& $\displaystyle n\left[
\left(\frac{\delta}{\det A}\right)^{1/n}
-\frac1n\log\frac{\delta}{\det A}-1
\right]$\\

$D_\beta$
& $\displaystyle
x_i\left(
\operatorname{arsinh}\frac{x_i}{\beta}
-\operatorname{arsinh}\frac{\lambda_i}{\beta}
\right)=\eta$
& $\displaystyle\sum_{i=1}^n\left[
f_\beta(x_i)-f_\beta(\lambda_i)
-(x_i-\lambda_i)
\operatorname{arsinh}\frac{\lambda_i}{\beta}
\right]$\\

$D_{\mathsf S}$
& $\displaystyle
\left(\frac{\delta}{\det A}\right)^{1/n}\lambda_i$
& $\displaystyle n\log
\frac{(\det A)^{1/n}+\delta^{1/n}}
{2(\delta\det A)^{1/(2n)}}$\\

$D_{\mathsf{QJS}}$
& $\displaystyle
\frac{x_i}{2}\log
\frac{2x_i}{\lambda_i+x_i}=\eta$
& $\displaystyle\frac12\sum_{i=1}^n\left[
\lambda_i\log\frac{2\lambda_i}{\lambda_i+x_i}
+x_i\log\frac{2x_i}{\lambda_i+x_i}
\right]$\\

$D_{\mathsf{BW}}$
& $\displaystyle
\frac{\bigl(\sqrt{\lambda_i}
\pm\sqrt{\lambda_i+4\eta}\bigr)^2}{4}$
& $\displaystyle
\sum_{i=1}^n(\sqrt{\lambda_i}-\sqrt{x_i})^2$\\

$D_{\mathsf{QH}}$
& $\displaystyle
\frac{\bigl(\sqrt{\lambda_i}
\pm\sqrt{\lambda_i+4\eta}\bigr)^2}{4}$
& $\displaystyle
\sum_{i=1}^n(\sqrt{\lambda_i}-\sqrt{x_i})^2$\\

$D_{\mathsf{Q},q}$
& $\displaystyle
\frac{q}{1-q}
\left(x_i-\lambda_i^{1-q}x_i^q\right)=\eta$
& $\displaystyle\frac1{1-q}\sum_{i=1}^n
\left[qx_i+(1-q)\lambda_i
-x_i^q\lambda_i^{1-q}\right]$\\

\end{tabu}
\end{adjustbox}
\caption{Solution to the left determinant-constrained nearness problem.}
\label{tab:det-left}
\end{table}
\end{proposition}

\begin{proof}
By Theorem~\ref{thm:reduce}, we reduce the problem to minimizing the reduced divergence over the domain
\[
\lambda(\mathcal S_\delta)
=\biggl\{x\in(0,\infty)^n:
x_1\geq\cdots\geq x_n,\,
\prod_{i=1}^n x_i=\delta\biggr\}.
\]
By \eqref{eq:mix}, we may drop the ordering constraint.

For every divergence under consideration, notice that $d_D(\lambda_i,x)$ and $d_D(x,\lambda_i)$ grow to $+\infty$ as $x\to\infty$. Since $\prod_i x_i=\delta$, every unbounded sequence of feasible eigenvalue vectors has an objective tending to $+\infty$. Consequently, the objective sublevel sets
are compact, and both minima are attained.

For $D_{\mathsf{vN}}$, forming the Lagrangian with the constraint $\log x_1+\dots+\log x_n=\log\delta$ yields 
\[
x_i-\lambda_i=\eta,\qquad
\prod_{i=1}^n(\lambda_i+\eta)=\delta.
\]
Since the product is strictly increasing for $\eta>-\lambda_n$, we obtain a unique solution.

For $D_{\mathsf{LD}}$, using the AM--GM inequality, we obtain
\[
\sum_{i=1}^n\frac{\lambda_i}{x_i}
\geq n\biggl(\frac{\det A}{\delta}\biggr)^{1/n},
\quad\text{and}\quad
\sum_{i=1}^n\frac{x_i}{\lambda_i}
\geq n\biggl(\frac{\delta}{\det A}\biggr)^{1/n},
\]
where equality is reached at $x_i=(\delta/\det A)^{1/n}\lambda_i$.

Similarly, for $D_{\mathsf S}$, Jensen's inequality applied to the strictly convex function $z\mapsto\log(1+e^z)$ gives
\[
\sum_{i=1}^n\log
\left(1+\frac{x_i}{\lambda_i}\right)
\geq n\log\left(1+
\left(\frac{\delta}{\det A}\right)^{1/n}\right),
\]
with equality at $x_i=(\delta/\det A)^{1/n}\lambda_i$.

For each remaining divergence, let $f_i(x)=d_D(\lambda_i,x)$ or $f_i(x)=d_D(x,\lambda_i)$. Forming the Lagrangian gives the equations
\begin{empheq}[left=\empheqlbrace]{align}
x_i f_i'(x_i)&=\eta,
\label{eq:det-stationarity}\\
\prod_{i=1}^n x_i&=\delta.
\label{eq:det-product}
\end{empheq}
The equations in Tables~\ref{tab:det-right} and~\ref{tab:det-left} are obtained by substituting the derivatives in Lemma~\ref{lem:separability}. We now show that their solutions determine a minimizer.

By Lemma~\ref{lem:separability}, if $\delta>\det A$, then $\eta>0$ and $x_i>\lambda_i$; if $\delta<\det A$, then $\eta<0$ and $x_i<\lambda_i$. If $\delta>\det A$, computing the first
derivative shows that $g_i(x)\coloneqq xf_i'(x)$ is strictly increasing on $(\lambda_i,\infty)$.
Therefore, for each $\eta>0$, \eqref{eq:det-stationarity} has a unique solution $x_i>\lambda_i$, which increases continuously with $\eta$. The constraint \eqref{eq:det-product} hence determines $\eta$ uniquely and gives a minimizer.

Suppose $\delta<\det A$. By direct calculations, we obtain that for $x\in(0,\lambda_i)$,
\[
g_i(x)<0,\qquad
\lim_{x\downarrow0}g_i(x)
=g_i(\lambda_i)=0,\qquad
g_i''(x)>0.
\]
Therefore, $g_i$ has a unique minimum on $(0,\lambda_i)$, and the equation $g_i(x)=\eta$ has at most two positive roots in this interval. Consider the change of variable $z_i=\log x_i$ for $i=1,\dots,n$. The Hessian of $F(z)=\sum_{i=1}^n f_i(e^{z_i})$ is
\[
\nabla^2F(z)
=\diag(x_1g_1'(x_1),\dots,
x_ng_n'(x_n)).
\]
At a local minimizer, $\nabla^2F$ is positive semidefinite on the tangent space
$\{z\in\R^n:\sum_i z_i=0\}$.

Whenever the two roots are distinct, $g_i'<0$ at the smaller root and $g_i'>0$ at the larger root. If two coordinates $i\neq j$ were distinct smaller roots, then $h=e_i-e_j$ would belong to the
tangent space and $h^\tp\nabla^2F(z)h
=x_i g_i'(x_i)+x_j g_j'(x_j)<0.$ Therefore, at most one coordinate can take a distinct smaller root.

Moreover, $\partial g_i(x)/\partial\lambda_i<0$. Implicit differentiation of
$g_i(x)=\eta$ yields
\[
\frac{\partial x}{\partial\lambda_i}
=-\frac{\partial_{\lambda_i}g_i(x)}
{g_i'(x)}.
\]
Thus, the smaller root decreases with $\lambda_i$, while the larger root increases. If a distinct smaller root occurred at $j<n$, then, since $\lambda_j\geq\lambda_n$ and the $n$th coordinate must take the larger or repeated root,
\[
x_j=x_-(\lambda_j)
\leq x_-(\lambda_n)
\leq x_+(\lambda_n)=x_n,
\]
with at least one strict inequality. This contradicts the nonincreasing ordering of the optimal solution. Therefore, a distinct smaller root can occur only at $i=n$.

Since a global minimizer exists and must satisfy \eqref{eq:det-stationarity} and \eqref{eq:det-product}, comparing the objective values of all admissible candidates with this root
structure yields a global minimizer. We recover the desired matrix by Corollary~\ref{cor:recover}.
\end{proof}

\subsection{Bosonic squeezed states}

Let $\p=\Sym(\C^n)$ and $\Omega=\Sym(\C^n)_{<1}$. For $0< k\leq n$ and $N>0$, we consider fidelity-optimal approximation subject to at most $k$ active squeezing modes and mean photon number at most $N$ \cite{HamiltonEtAl2017,KruseEtAl2019}:
\begin{equation}\label{eq:boson-budget}
\min_{X\in\mathcal S_{k,N}}C_{\mathsf B}(A,X),\qquad
\mathcal S_{k,N}\coloneqq
\bigl\{X\in\Omega:\rank(X)\leq k,\quad
\tr\bigl[X^*X(I-X^*X)^{-1}\bigr]\leq N\bigr\}.
\end{equation}

\begin{proposition}\label{prop:boson-budget}
For $A \in \Sym(\C^n)_{<1}$, let $A=U\diag(\sigma_1,\dots,\sigma_n)U^\tp$ be its Autonne--Takagi
decomposition. If $N\geq \sum_{i=1}^k \sigma_i^2/(1-\sigma_i^2)$, let $v_i=\sigma_i$ for $i=1,\dots, k$ and $v_i =0$ for $i = k+1,\dots,n$. If $0<N<\sum_{i=1}^k \sigma_i^2/(1-\sigma_i^2)$, let $\eta>0$ be the unique parameter such that for $i=1,\dots,k$, $v_i\in[0,\sigma_i]$ is the unique root of
\begin{equation}\label{eq:boson-cubic}
x^3+(2\eta-1)\sigma_ix^2-(1+2\eta)x+\sigma_i=0.
\end{equation}
and satisfies
\begin{equation}\label{eq:boson-eta}
\sum_{i=1}^k\frac{v_i^2}{1-v_i^2}=N.
\end{equation}
Then $X_*\coloneqq U\diag(v_1,\dots,v_n)U^\tp$
is a minimizer of \eqref{eq:boson-budget}, with 
\[
\min_{X\in\mathcal S_{k,N}}C_{\mathsf B}(A,X) = 
\begin{cases}
    \sum_{i=1}^n\log\frac{1-\sigma_iv_i}{\sqrt{(1-\sigma_i^2)(1-v_i^2)}} & \text{if } 0<N<\sum_{i=1}^k\frac{\sigma_i^2}{1-\sigma_i^2},\\
    -\frac12\sum_{i=k+1}^n\log(1-\sigma_i^2) & \text{if } N\geq\sum_{i=1}^k\frac{\sigma_i^2}{1-\sigma_i^2}.
\end{cases}
\]
\end{proposition}

\begin{proof}
Since $\mathcal S_{k,N}$ is $\U(n)$-invariant, by Theorem~\ref{thm:reduce}, \eqref{eq:boson-budget} reduces to 
\[
\min\biggl\{\widehat C_{\mathsf B}(\sigma,v):
1>v_1\geq\cdots\geq v_k\geq
v_{k+1}=\cdots=v_n=0,\,
\sum_{i=1}^k\frac{v_i^2}{1-v_i^2}\leq N\biggr\}.
\]
Consider the change of variable $a_i=\operatorname{artanh}\sigma_i$ and $t_i=\operatorname{artanh}v_i$ for $i\leq k$. Then, removing the constant terms, we solve the relaxed problem 
\begin{equation}\label{eq:boson-convex}
    \min_{t \in \R^k}\,\sum_{i=1}^k\log\cosh(a_i-t_i)\quad \text{such that} \quad \sum_{i=1}^k\sinh^2t_i\leq N.
\end{equation}
Since
\[
\frac{d^2}{dt^2}\log\cosh(a_i-t)=\operatorname{sech}^2(a_i-t)>0,
\qquad
\frac{d^2}{dt^2}\sinh^2t=2\cosh(2t)>0,
\]
\eqref{eq:boson-convex} is a strictly convex optimization problem. We now show that its minimizer automatically satisfies the omitted ordering and nonnegativity constraints.

Since $t=0$ satisfies Slater's condition, the KKT conditions are necessary and sufficient for optimality. Forming the Lagrangian
\[
\sum_{i=1}^k\log\cosh(a_i-t_i)
+\eta\biggl(\sum_{i=1}^k\sinh^2t_i-N\biggr),
\]
where $\eta\geq0$, the KKT conditions are
\begin{align}
    \tanh(t_i-a_i)+\eta\sinh(2t_i)&=0,
\qquad i=1,\dots,k, \label{eq:stationarity}\\
\eta\biggl(\sum_{i=1}^k\sinh^2t_i-N\biggr)&=0,
\qquad \sum_{i=1}^k\sinh^2t_i\leq N \notag.
\end{align}
Notice that for fixed $\eta\geq0$, the function $t\mapsto\tanh(t-a_i)+\eta\sinh(2t)$ is strictly increasing by computing its first derivative. Moreover, its values at $t=0$ and $t=a_i$ are $-\tanh a_i\leq0$ and $\eta\sinh(2a_i)\geq0$, respectively. Hence, there exists an unique $t_i$ such that $\tanh(t_i-a_i)+\eta\sinh(2t_i)=0$ for each $i=1,\dots,k$. Furthermore, since the $\tanh(t-a_i)+\eta\sinh(2t)$ is decreasing in $a_i$, the solution $t_i$ is nondecreasing in $a_i$. Therefore,$t_1\geq\cdots\geq t_k\geq0$.

If $N\geq\sum_{i=1}^k\sinh^2a_i$, then $\eta=0$ and $t_i=a_i$ suffice. Otherwise, $0<N<\sum_{i=1}^k\sinh^2a_i$, so $\eta>0$ and $\sum_{i=1}^k\sinh^2t_i=N$ by complementary slackness. For fixed $a_i$, implicitly differentiating \eqref{eq:stationarity} with respect to $\eta$ yields
\[
\frac{dt_i}{d\eta}
=-\frac{\sinh(2t_i)}
{\operatorname{sech}^2(t_i-a_i)+2\eta\cosh(2t_i)}
\leq0,
\]
with strict inequality whenever $a_i>0$. Moreover, since $t_i=a_i$ for $\eta = 0$, and $0\leq\eta\sinh(2t_i) =\tanh(a_i-t_i)\leq1$, $t_i\to0$ as $\eta\to+\infty$. Consequently, $\eta\mapsto\sum_{i=1}^k\sinh^2t_i$, where $t_i$ depends on $\eta$, is continuous and strictly decreasing from
$\sum_{i=1}^k\sinh^2a_i$ to $0$. Therefore, there exists a unique $\eta>0$ satisfying $\sum_{i=1}^k\sinh^2t_i=N$.

Finally, substituting $v_i=\tanh t_i$ and $\sigma_i=\tanh a_i$ into \eqref{eq:stationarity} gives \eqref{eq:boson-cubic}. The minimum value follows by substitution into the reduced divergence, and $X_*$ is recovered by Corollary~\ref{cor:recover}.
\end{proof}

\subsection{Fermionic pairing states}

Let $\p=\Alt(\C^{2n})$ and $\Omega=\p_{<1}$. For $0< k\leq n$ and $N> 0$, we consider fidelity-optimal approximation subject to at most $k$ active pairing channels and mean fermion number at most $N$ \cite{BlochMessiah1962,TsuchimochiScuseria2009}:
\begin{equation}\label{eq:fermion-budget}
\min_{X\in\mathcal S_{k,N}}C_{\mathsf F}(A,X),\qquad
\mathcal S_{k,N}\coloneqq
\bigl\{X\in\Omega:\rank(X)\leq2k,\quad
\tr\bigl[X^*X(I+X^*X)^{-1}\bigr]\leq N\bigr\}.
\end{equation}

\begin{proposition}\label{prop:fermion-budget}
Let $A=Q(\diag(\sigma_1,\dots,\sigma_n)\otimes J)Q^\tp \neq 0$ be its Youla decomposition, where $1>\sigma_1\geq\cdots\geq\sigma_n\geq0$. Let $v_i=0$ for $i>k$. If $N\geq2\sum_{i=1}^k\sigma_i^2/(1+\sigma_i^2)$, let $v_i=\sigma_i$ for $i=1,\dots,k$.
Otherwise, for
$0<N<2\sum_{i=1}^k\sigma_i^2/(1+\sigma_i^2)$,
let $\eta>0$ be the unique parameter such that
for $i=1,\dots,k$, $v_i\in[0,\sigma_i]$ is the unique root of
\begin{equation}\label{eq:fermion-cubic}
x^3+(2\eta-1)\sigma_i x^2+(1+2\eta)x-\sigma_i=0
\end{equation}
and satisfies
\begin{equation}\label{eq:fermion-eta}
2\sum_{i=1}^k\frac{v_i^2}{1+v_i^2}=N.
\end{equation}
Then
\[
X_*\coloneqq Q\bigl(\diag(v_1,\dots,v_n)\otimes J\bigr)Q^\tp
\]
is a minimizer of \eqref{eq:fermion-budget}, with
\[
\min_{X\in\mathcal S_{k,N}}C_{\mathsf F}(A,X)
=
\begin{cases}
\displaystyle
\sum_{i=1}^n
\log\frac{(1+\sigma_i^2)(1+v_i^2)}
{(1+\sigma_iv_i)^2},
&\displaystyle
0 < N<2\sum_{i=1}^k\frac{\sigma_i^2}{1+\sigma_i^2},
\\[8pt]
\displaystyle
\sum_{i=k+1}^n\log(1+\sigma_i^2),
&\displaystyle
N\geq2\sum_{i=1}^k\frac{\sigma_i^2}{1+\sigma_i^2}.
\end{cases}
\]
\end{proposition}

\begin{proof}
Since $\mathcal S_{k,N}$ is $\U(2n)$-invariant, by Theorem~\ref{thm:reduce}, \eqref{eq:fermion-budget} reduces to
\[
\min\biggl\{\widehat C_{\mathsf F}(\sigma,v):
v\in\rho(\Omega),\;
v_i=0\ (i>k),\;
2\sum_{i=1}^k\frac{v_i^2}{1+v_i^2}\leq N\biggr\}.
\]
Consider the change of variables
$a_i=\arctan\sigma_i$ and $t_i=\arctan v_i$
for $i=1,\dots,k$.
Then, removing the constant terms and temporarily omitting the ordering and nonnegativity constraints, we consider the relaxed problem
\begin{equation}\label{eq:fermion-convex}
\min_{t\in[-\pi/4,\pi/4]^k}
-2\sum_{i=1}^k\log\cos(a_i-t_i)
\quad\text{such that}\quad
2\sum_{i=1}^k\sin^2t_i\leq N.
\end{equation}
Since
\[
\frac{d^2}{dt^2}\bigl[-2\log\cos(a_i-t)\bigr]
=2\sec^2(a_i-t)>0,
\qquad
\frac{d^2}{dt^2}\bigl[2\sin^2t\bigr]
=4\cos(2t)\geq0
\]
for $t\in[-\pi/4,\pi/4]$, \eqref{eq:fermion-convex} is a strictly convex optimization problem. Since $N>0$, Slater's condition holds For an interior point, the stationarity condition is
\begin{equation}\label{eq:stationarity2}
    \tan(t_i-a_i)+\eta\sin(2t_i)=0,
\qquad i=1,\dots,k,
\end{equation}
where $\eta\geq0$, together with
\[
\eta\biggl(2\sum_{i=1}^k\sin^2t_i-N\biggr)=0.
\]
For fixed $\eta\geq0$, the left-hand side of \eqref{eq:stationarity2} is strictly increasing on $[0,a_i]$ as $\sec^2(t-a_i)+2\eta\cos(2t)>0$. Its values at $t=0$ and $t=a_i$ have opposite weak signs, so there exists a unique solution $t_i\in[0,a_i]$. Moreover, implicitly differentiating \eqref{eq:stationarity2} yields that
\[
\begin{aligned}
\frac{dt_i}{da_i}
&=\frac{\sec^2(t_i-a_i)}
{\sec^2(t_i-a_i)+2\eta\cos(2t_i)}>0,\\
\frac{dt_i}{d\eta}
&=-\frac{\sin(2t_i)}
{\sec^2(t_i-a_i)+2\eta\cos(2t_i)}\leq0.
\end{aligned}
\]
Thus $t_1\geq\cdots\geq t_k\geq0$, so the minimizer satisfies the omitted constraints. Since $t_i\leq a_i<\pi/4$, these solutions lie in $[-\pi/4,\pi/4]^k$. Furthermore, since $0\leq\eta\sin(2t_i)=\tan(a_i-t_i)\leq\tan a_i$, $t_i\to0$ as $\eta\to+\infty$. As in the proof of Proposition~\ref{prop:boson-budget}, the particle-number constraint therefore determines a unique $\eta>0$ when
$0<N<2\sum_{i=1}^k\sin^2a_i$; otherwise $\eta=0$ and $t_i=a_i$.

Finally, substituting $v_i=\tan t_i$ and $\sigma_i=\tan a_i$ into the stationarity equation gives \eqref{eq:fermion-cubic}. The minimum value follows by substitution into the reduced divergence, and $X_*$ is recovered by Corollary~\ref{cor:recover}.
\end{proof}

\section{Conclusion}

To answer the question posed at the beginning of the article, we proposed the notion of isotropic divergences, under which a solution to a matrix nearness problem can be obtained from merely manipulating the eigenvalues and singular values of the approximand whenever the constraints respect the Lie group symmetry. We then verify the isotropy of a range of divergences arising from applications. As illustrated by the examples in Section~\ref{sec:solve}, we hope this work provides the tools to solve future nearness problems systematically.

\subsection*{AI disclosure} This work is initiated with the author's work on \cite{ten} where the author gradually formed the definition of isotropic divergences. In the process of explorations, ChatGPT models released in 2025 and 2026 have been used for the purpose of verifying proofs, checking calculations, and literature search, leading to this article. The author welcomes questions regarding the role of AI in this work, as well as comments and suggestions concerning best practices for the disclosure of AI-assisted mathematical research.


\begin{thebibliography}{10}

\bibitem{Abe2003}
S.~Abe.
\newblock Nonadditive generalization of the quantum kullback-leibler divergence for measuring the degree of purification.
\newblock {\em Phys. Rev. A}, 68:032302, Sep 2003.

\bibitem{Amari2016}
S.-i. Amari.
\newblock {\em Information geometry and its applications}, volume 194 of {\em Applied Mathematical Sciences}.
\newblock Springer, [Tokyo], 2016.

\bibitem{Anderson2003}
T.~W. Anderson.
\newblock {\em An introduction to multivariate statistical analysis}.
\newblock Wiley Series in Probability and Statistics. Wiley-Interscience [John Wiley \& Sons], Hoboken, NJ, third edition, 2003.

\bibitem{Atiyah}
M.~F. Atiyah.
\newblock Convexity and commuting {H}amiltonians.
\newblock {\em Bull. London Math. Soc.}, 14(1):1--15, 1982.

\bibitem{Autonne15}
L.~Autonne.
\newblock Sur les matrices hypohermitiennes et sur les matrices unitaires.
\newblock {\em Ann. Univ. Lyon Nouv. S{\'e}r. I}, 38:1--77, 1915.

\bibitem{BachJordan2005}
F.~R. Bach and M.~I. Jordan.
\newblock A probabilistic interpretation of canonical correlation analysis.
\newblock Technical Report 688, Department of Statistics, University of California, Berkeley, 2005.

\bibitem{BhatiaGaubertJain2019}
R.~Bhatia, S.~Gaubert, and T.~Jain.
\newblock Matrix versions of the {H}ellinger distance.
\newblock {\em Lett. Math. Phys.}, 109(8):1777--1804, 2019.

\bibitem{BhatiaJainLim2019}
R.~Bhatia, T.~Jain, and Y.~Lim.
\newblock On the {B}ures-{W}asserstein distance between positive definite matrices.
\newblock {\em Expo. Math.}, 37(2):165--191, 2019.

\bibitem{BianchiHacklYokomizo2015}
E.~Bianchi, L.~Hackl, and N.~Yokomizo.
\newblock Entanglement entropy of squeezed vacua on a lattice.
\newblock {\em Phys. Rev. D}, 92(8):085045, 20, 2015.

\bibitem{BlochMessiah1962}
C.~Bloch and A.~Messiah.
\newblock The canonical form of an antisymmetric tensor and its application to the theory of superconductivity.
\newblock {\em Nuclear Phys.}, 39:95--106, 1962.

\bibitem{Calabi1953}
E.~Calabi.
\newblock Isometric imbedding of complex manifolds.
\newblock {\em Ann. of Math. (2)}, 58:1--23, 1953.

\bibitem{CarlssonRotureau2021}
B.~G. Carlsson and J.~Rotureau.
\newblock New and practical formulation for overlaps of bogoliubov vacua.
\newblock {\em Phys. Rev. Lett.}, 126:172501, Apr 2021.

\bibitem{ChenPedersenAndersen2025}
L.~Chen, F.~H. Pedersen, and M.~S. Andersen.
\newblock Matrix nearness problems with off-block-diagonal rank constraints.
\newblock {\em Linear Algebra Appl.}, 721:635--673, 2025.

\bibitem{Cherian2013}
A.~Cherian, S.~Sra, A.~Banerjee, and N.~Papanikolopoulos.
\newblock Jensen-{B}regman logdet divergence with application to efficient similarity search for covariance matrices.
\newblock {\em IEEE Transactions on Pattern Analysis and Machine Intelligence}, 35(9):2161--2174, 2013.

\bibitem{logdet2}
J.~V. Davis, B.~Kulis, P.~Jain, S.~Sra, and I.~S. Dhillon.
\newblock Information-theoretic metric learning.
\newblock In {\em Proceedings of the 24th International Conference on Machine Learning}, ICML '07, pages 209--216, New York, NY, USA, 2007. Association for Computing Machinery.

\bibitem{logdet1}
D.~K. Dey and C.~Srinivasan.
\newblock Estimation of a covariance matrix under {S}tein's loss.
\newblock {\em Ann. Statist.}, 13(4):1581--1591, 1985.

\bibitem{DT2007}
I.~S. Dhillon and J.~A. Tropp.
\newblock Matrix nearness problems with {B}regman divergences.
\newblock {\em SIAM J. Matrix Anal. Appl.}, 29(4):1120--1146, 2007.

\bibitem{Duchi2008}
J.~C. Duchi, S.~Shalev-Shwartz, Y.~Singer, and T.~D. Chandra.
\newblock Efficient projections onto the l1-ball for learning in high dimensions.
\newblock In {\em International Conference on Machine Learning}, 2008.

\bibitem{Eck}
C.~Eckart and G.~Young.
\newblock The approximation of one matrix by another of lower rank.
\newblock {\em Psychometrika}, 1(3):211--218, 1936.

\bibitem{FanHoffman1955}
K.~Fan and A.~J. Hoffman.
\newblock Some metric inequalities in the space of matrices.
\newblock {\em Proc. Amer. Math. Soc.}, 6:111--116, 1955.

\bibitem{FeliceAy2019}
D.~Felice and N.~Ay.
\newblock Canonical divergence for flat $\alpha$-connections: Classical and quantum.
\newblock {\em Entropy}, 21(9), 2019.

\bibitem{Frisch1929}
R.~Frisch.
\newblock Correlation and scatter in statistical variables.
\newblock {\em Nordic Statistical Journal}, 1:36--102, 1929.

\bibitem{Determinantal}
A.~Fr\"uhbis-Kr\"uger and M.~Zach.
\newblock Determinantal singularities.
\newblock In {\em Handbook of geometry and topology of singularities {IV}}, pages 45--159. Springer, Cham, [2023] \copyright 2023.

\bibitem{FuruichiYanagiKuriyama2004}
S.~Furuichi, K.~Yanagi, and K.~Kuriyama.
\newblock Fundamental properties of {T}sallis relative entropy.
\newblock {\em J. Math. Phys.}, 45(12):4868--4877, 2004.

\bibitem{hypentropy}
U.~Ghai, E.~Hazan, and Y.~Singer.
\newblock Exponentiated gradient meets gradient descent.
\newblock In A.~Kontorovich and G.~Neu, editors, {\em Proceedings of the 31st International Conference on Algorithmic Learning Theory}, volume 117 of {\em Proceedings of Machine Learning Research}, pages 386--407. PMLR, 08 Feb--11 Feb 2020.

\bibitem{GVL1980}
G.~H. Golub and C.~F. Van~Loan.
\newblock An analysis of the total least squares problem.
\newblock {\em SIAM J. Numer. Anal.}, 17(6):883--893, 1980.

\bibitem{HamiltonEtAl2017}
C.~S. Hamilton, R.~Kruse, L.~Sansoni, S.~Barkhofen, C.~Silberhorn, and I.~Jex.
\newblock Gaussian boson sampling.
\newblock {\em Phys. Rev. Lett.}, 119:170501, Oct 2017.

\bibitem{Helgason2001}
S.~Helgason.
\newblock {\em Differential geometry, {L}ie groups, and symmetric spaces}, volume~34 of {\em Graduate Studies in Mathematics}.
\newblock American Mathematical Society, Providence, RI, 2001.
\newblock Corrected reprint of the 1978 original.

\bibitem{HiaiEtAl2011}
F.~Hiai, M.~Mosonyi, D.~Petz, and C.~B\'eny.
\newblock Quantum {$f$}-divergences and error correction.
\newblock {\em Rev. Math. Phys.}, 23(7):691--747, 2011.

\bibitem{Higham1}
N.~J. Higham.
\newblock Computing a nearest symmetric positive semidefinite matrix.
\newblock {\em Linear Algebra Appl.}, 103:103--118, 1988.

\bibitem{Higham2}
N.~J. Higham.
\newblock Matrix nearness problems and applications.
\newblock In {\em Applications of matrix theory ({B}radford, 1988)}, volume~22 of {\em Inst. Math. Appl. Conf. Ser. New Ser.}, pages 1--27. Oxford Univ. Press, New York, 1989.

\bibitem{Horn}
A.~Horn.
\newblock Doubly stochastic matrices and the diagonal of a rotation matrix.
\newblock {\em Amer. J. Math.}, 76:620--630, 1954.

\bibitem{HJ1991}
R.~A. Horn and C.~R. Johnson.
\newblock {\em Topics in matrix analysis}.
\newblock Cambridge University Press, Cambridge, 1991.

\bibitem{HZ1012}
R.~A. Horn and F.~Zhang.
\newblock A generalization of the complex {A}utonne-{T}akagi factorization to quaternion matrices.
\newblock {\em Linear Multilinear Algebra}, 60(11-12):1239--1244, 2012.

\bibitem{Hotelling1936}
H.~Hotelling.
\newblock Relations between two sets of variates.
\newblock {\em Biometrika}, 28(3/4):321--377, 1936.

\bibitem{JS2026}
P.~Jaap and O.~Sander.
\newblock How to project onto {${\rm SL}(n)$}.
\newblock {\em J. Optim. Theory Appl.}, 209(2):Paper No. 40, 43, 2026.

\bibitem{Knapp2002}
A.~W. Knapp.
\newblock {\em Lie groups beyond an introduction}, volume 140 of {\em Progress in Mathematics}.
\newblock Birkh\"auser Boston, Inc., Boston, MA, second edition, 2002.

\bibitem{Kostant}
B.~Kostant.
\newblock On convexity, the {W}eyl group and the {I}wasawa decomposition.
\newblock {\em Ann. Sci. \'Ecole Norm. Sup. (4)}, 6:413--455, 1973.

\bibitem{Krichene2015}
W.~Krichene, S.~Krichene, and A.~Bayen.
\newblock Efficient bregman projections onto the simplex.
\newblock In {\em 2015 54th IEEE Conference on Decision and Control (CDC)}, pages 3291--3298, 2015.

\bibitem{KruseEtAl2019}
R.~Kruse, C.~S. Hamilton, L.~Sansoni, S.~Barkhofen, C.~Silberhorn, and I.~Jex.
\newblock Detailed study of gaussian boson sampling.
\newblock {\em Phys. Rev. A}, 100:032326, Sep 2019.

\bibitem{KL}
S.~Kullback and R.~A. Leibler.
\newblock On information and sufficiency.
\newblock {\em Ann. Math. Statistics}, 22:79--86, 1951.

\bibitem{Lewis1996}
A.~S. Lewis.
\newblock Group invariance and convex matrix analysis.
\newblock {\em SIAM J. Matrix Anal. Appl.}, 17(4):927--949, 1996.

\bibitem{lewis2000convex}
A.~S. Lewis.
\newblock Convex analysis on {C}artan subspaces.
\newblock {\em Nonlinear Anal.}, 42(5):813--820, 2000.

\bibitem{LuoZhang2004}
S.~Luo and Q.~Zhang.
\newblock Informational distance on quantum-state space.
\newblock {\em Phys. Rev. A (3)}, 69(3):032106, 8, 2004.

\bibitem{Majtey2005}
A.~P. Majtey, P.~W. Lamberti, and D.~P. Prato.
\newblock Jensen-{S}hannon divergence as a measure of distinguishability between mixed quantum states.
\newblock {\em Phys. Rev. A}, 72:052310, Nov 2005.

\bibitem{mir}
L.~Mirsky.
\newblock Symmetric gauge functions and unitarily invariant norms.
\newblock {\em Quart. J. Math. Oxford Ser. (2)}, 11:50--59, 1960.

\bibitem{MW1931}
F.~D. Murnaghan and A.~Wintner.
\newblock A canonical form for real matrices under orthogonal transformations.
\newblock {\em Proceedings of the National Academy of Sciences of the United States of America}, 17(7):417--420, 1931.

\bibitem{NielsenBoltz2011}
F.~Nielsen and S.~Boltz.
\newblock The burbea-rao and bhattacharyya centroids.
\newblock {\em IEEE Trans. Inf. Theor.}, 57(8):5455--5466, Aug. 2011.

\bibitem{NielsonChuang2000}
M.~A. Nielsen and I.~L. Chuang.
\newblock {\em Quantum computation and quantum information}.
\newblock Cambridge University Press, Cambridge, 2000.

\bibitem{Nievergelt}
Y.~Nievergelt.
\newblock Schmidt-{M}irsky matrix approximation with linearly constrained singular values.
\newblock {\em Linear Algebra Appl.}, 261:207--219, 1997.

\bibitem{Robledo2009}
L.~M. Robledo.
\newblock Sign of the overlap of hartree-fock-bogoliubov wave functions.
\newblock {\em Phys. Rev. C}, 79:021302(R), Feb 2009.

\bibitem{Schur}
I.~Schur.
\newblock {\"U}ber eine klasse von mittelbildungen mit anwendungen auf die determinantentheorie.
\newblock {\em Sitzungsber. Berl. Math. Ges.}, 22:9--20, 1923.

\bibitem{Smolin2012}
J.~A. Smolin, J.~M. Gambetta, and G.~Smith.
\newblock Efficient method for computing the maximum-likelihood quantum state from measurements with additive gaussian noise.
\newblock {\em Phys. Rev. Lett.}, 108:070502, Feb 2012.

\bibitem{Sra2016}
S.~Sra.
\newblock Positive definite matrices and the {S}-divergence.
\newblock {\em Proc. Amer. Math. Soc.}, 144(7):2787--2797, 2016.

\bibitem{Sra2021}
S.~Sra.
\newblock Metrics induced by {J}ensen-{S}hannon and related divergences on positive definite matrices.
\newblock {\em Linear Algebra Appl.}, 616:125--138, 2021.

\bibitem{Takagi1924}
T.~Takagi.
\newblock On an algebraic problem reluted to an analytic theorem of carath{\'e}odory and fej{\'e}r and on an allied theorem of landau.
\newblock {\em Japanese journal of mathematics :transactions and abstracts}, 1:83--93, 1924.

\bibitem{Tam1999}
T.-Y. Tam.
\newblock An extension of a result of {L}ewis.
\newblock {\em Electron. J. Linear Algebra}, 5:1--10, 1999.

\bibitem{TsuchimochiScuseria2009}
T.~Tsuchimochi and G.~E. Scuseria.
\newblock Strong correlations via constrained-pairing mean-field theory.
\newblock {\em The Journal of Chemical Physics}, 131(12):121102, 09 2009.

\bibitem{Wald1943}
A.~Wald.
\newblock On the efficient design of statistical investigations.
\newblock {\em Ann. Math. Statistics}, 14:134--140, 1943.

\bibitem{ten}
R.~T. Wang, X.~Lu, K.~Ye, and L.-H. Lim.
\newblock Singular value decomposition on manifolds.
\newblock In preparation.

\bibitem{Wilkinson1986}
J.~H. Wilkinson.
\newblock Sensitivity of eigenvalues. {II}.
\newblock {\em Utilitas Math.}, 30:243--286, 1986.

\bibitem{Wilks1932}
S.~S. Wilks.
\newblock Certain generalizations in the analysis of variance.
\newblock {\em Biometrika}, 24(3/4):471--494, 1932.

\bibitem{Youla}
D.~C. Youla.
\newblock A normal form for a matrix under the unitary congruence group.
\newblock {\em Canadian Journal of Mathematics}, 13:694--704, 1961.

\end{thebibliography}
\end{document}